\documentclass[a4paper,11pt]{amsart}
\usepackage{amssymb,mathpazo,bm}
\usepackage[shortlabels]{enumitem} 
\usepackage[hmargin=2.5cm,vmargin=2cm]{geometry}
\usepackage[all,poly,arc]{xy}

\newcommand{\ppq}{\leqslant}
\newcommand{\pgq}{\geqslant}

\newcommand{\dsty}{\displaystyle}
\newcommand{\tsty}{\textstyle}

\newcommand{\End}{\mathrm{End}}

\newcommand{\im}{\operatorname{Im}\nolimits}
\renewcommand{\ker}{\operatorname{Ker}\nolimits}
\newcommand{\car}{\operatorname{char}\nolimits}
\newcommand{\Hom}{\operatorname{Hom}\nolimits}
\newcommand{\Ext}{\operatorname{Ext}\nolimits}
\newcommand{\HH}{\operatorname{HH}\nolimits}

\newcommand{\rad}{\operatorname{rad}\nolimits}
\newcommand{\soc}{\operatorname{soc}\nolimits}
\renewcommand{\top}{\operatorname{top}\nolimits}

\newcommand{\ot}{\otimes}

\newcommand{\mo}{\mathfrak{o}}
\newcommand{\mt}{\mathfrak{t}}
\newcommand{\rrad}{\mathfrak{r}}

\newcommand{\K}{\mathcal{K}}

\newcommand{\B}{\mathcal{B}}

\newcommand{\D}{\mathcal{D}}

\newcommand{\R}{\mathcal{R}}

\newcommand{\Q}{\mathcal{Q}}

\newcommand{\A}{\mathcal{A}}

\newcommand{\set}[1]{\left\{ #1 \right\}}

\newtheorem{theorem}{{Theorem}}[section]
\newenvironment{thm}{\begin{theorem}}{\end{theorem}}
\newcommand{\bt}{\begin{thm}}
\newcommand{\et}{\end{thm}}

\newtheorem{corollaire}[theorem]{{Corollary}}
\newenvironment{cor}{\begin{corollaire}}{\end{corollaire}}
\newcommand{\bc}{\begin{cor}}
\newcommand{\ec}{\end{cor}}

\newtheorem{lemme}[theorem]{{Lemma}}
\newenvironment{lemma}{\begin{lemme}}{\end{lemme}}
\newcommand{\bl}{\begin{lemma}}
\newcommand{\el}{\end{lemma}}

\newtheorem{proposition}[theorem]{{Proposition}}
\newcommand{\bprop}{\begin{proposition}}
\newcommand{\eprop}{\end{proposition}}

\newtheorem{definition}[theorem]{{Definition}}
\newenvironment{dfn}{\begin{definition} \rm}{\end{definition}}
\newcommand{\bd}{\begin{dfn}}
\newcommand{\ed}{\end{dfn}}

\newtheorem*{remark}{Remark}
\newcommand{\br}{\begin{remark}}
\newcommand{\er}{\end{remark}}

\begin{document}

\title[Classification of symmetric special biserial algebras with \dots]
{Classification of symmetric special biserial
algebras with at most one non-uniserial indecomposable projective}
\author[Snashall]{Nicole Snashall}
\address{Nicole Snashall\\Department of Mathematics\\ University of Leicester\\
University Road\\ Leicester, LE1 7RH\\ England}
\email{N.Snashall@mcs.le.ac.uk}
\author[Taillefer]{Rachel Taillefer}
\address{Rachel Taillefer\\ Clermont Universit\'e, Universit\'e Blaise Pascal, Laboratoire de Math\'ematiques\\ BP 10448, F-63000 Clermont-Ferrand -
CNRS, UMR 6620, Laboratoire de Math\'ematiques, F-63177 Aubi\`ere}
\email{Rachel.Taillefer@math.univ-bpclermont.fr}


\begin{abstract}
We consider a natural generalisation of symmetric Nakayama algebras, namely, symmetric special
biserial algebras with at most one non-uniserial indecomposable projective module. We describe the
basic algebras explicitly by quiver and relations, then classify them up to derived equivalence and
up to stable equivalence of Morita type. This includes the algebras of \cite{BHS}, where they study
the weakly symmetric algebras of Euclidean type, as well as some algebras of dihedral type.
\end{abstract}

\date{\today} \maketitle

\section*{Introduction}

Let $K$ be an algebraically closed field. We consider in this paper a generalisation of symmetric Nakayama
$K$-algebras.
A symmetric Nakayama $K$-algebra is a symmetric $K$-algebra $A$ such that all indecomposable projective modules are
uniserial. These algebras are well-known and have been classified up to Morita equivalence: every symmetric Nakayama
algebra is Morita equivalent to exactly one algebra $N_m^n$ defined by the quiver
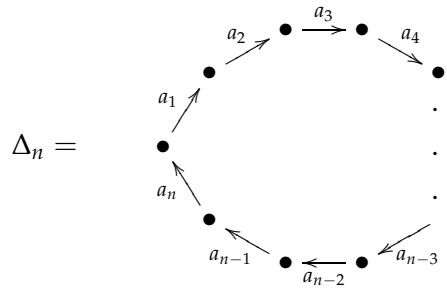
\[\Delta_n=\vcenter{\xymatrix@C=.4pc@R=.4pc{
  && &&& \bullet \ar^{a_3}[rr] && \bullet \ar^{a_4}[rrd] \\
  &
       && \bullet \ar^{a_2}[urr] &&&&&& \bullet
  \ar@{*{}*{\rule[2.5mm]{0mm}{0mm}.\rule[-2.5mm]{0mm}{0mm}}*{}}[dddd]
      \\
   &&&&&&&&&&
    \\
   && {\bullet} \ar^{a_1}[ruu]
      \\
    &&&&&&&&&&
     \\
   &
         && \bullet \ar^{a_{n}}[luu]
        &&&&&& \ar^{a_{n-3}}[lld] \\
  && &&&
      \bullet \ar^{a_{n-1}}[llu] && \bullet \ar^{a_{n-2}}[ll] \\
}}\] and the ideal of relations $L_m$ in $K\Delta_n$ generated by all paths of length $nm+1.$ Note that in particular,
the basic algebra associated to $A$ is special biserial.

Our aim is to describe the basic indecomposable finite-dimensional $K$-algebras $A$ which are symmetric special
biserial algebras with at most one non-uniserial indecomposable projective module. These algebras
include the symmetric Nakayama algebras,  certain algebras in \cite{BHS} that occur in the
classification, up to derived equivalence, of all weakly symmetric algebras of Euclidean type,
as well as some algebras of dihedral type, see \cite{E}. In
this paper we also distinguish, up to derived equivalence and up to stable equivalence of Morita
type, the basic indecomposable finite-dimensional symmetric special biserial algebras which have at
most one non-uniserial indecomposable projective module. It is well-known that all special biserial
algebras are tame \cite{WW}. Morever, it was proved by Al-Nofayee in \cite{AlN} (and by Rickard \cite{Ri} for
the symmetric case) that if $A$ and $B$ are derived equivalent algebras, then $A$ is selfinjective
if and only if $B$ is selfinjective. It was also proved by Pogorza\l y \cite{Po} that if $A$ is a selfinjective
special biserial algebra that is not a Nakayama algebra and if $A$ and $B$ are stably equivalent of Morita
type, then
$B$ is also a selfinjective
special biserial algebra. The algebras given in \cite{BHS} are Brauer graph algebras, and we recall that Brauer tree algebras play an important role in the Morita equivalence classification of blocks of group algebras of finite type (see \cite{A,B}). We use the theory of generalised Brauer tree algebras as part of the classification of our algebras up to derived equivalence. We refer the reader also to \cite{Sk}, where Skowro\'nski discusses the extensive programme to determine the derived equivalence classes of all tame selfinjective algebras.

\bigskip

We begin this paper with some background and properties about basic symmetric algebras, so that, in
Section~\ref{sec:quiver}, we can describe by quiver and relations all basic indecomposable
finite-dimensional algebras which are symmetric special biserial algebras with at most one
non-uniserial indecomposable projective module. In order to distinguish our algebras up to derived
equivalence and up to stable equivalence of Morita type, we use several invariants including
Hochschild cohomology which we discuss in  Section~\ref{sec:HH}. Section~\ref{sec:final} contains
the full classification of our algebras up to derived equivalence, and in addition to the dimensions
of the Hochschild cohomology groups, we use Cartan invariants (see \cite[Proposition 1.5]{BS} for a
proof of derived invariance) and K\"ulshammer invariants (or generalised Reynolds ideals, whose
derived invariance was proved in \cite{Z}). The final section gives the full classification of our
algebras up to stable equivalence of Morita type, based on the classification up to derived
equivalence of Section~\ref{sec:final} and using similar invariants.

We assume throughout that $A$ is a basic indecomposable finite-dimensional algebra over the algebraically closed field
$K$ so that $A$ is isomorphic to $K{\mathcal Q}/I$ for some unique connected quiver ${\mathcal Q}$ and admissible ideal $I$ of $K{\mathcal Q}$. We let $\rad(A)$ denote the Jacobson radical of $A$.

For any two positive integers $p$ and $q$ with $p\ppq q$, we define the quiver $\Q_{(p,q)}$ to be the quiver formed of
two oriented cycles, of lengths $p$ and $q$ respectively, joined at one vertex labelled $1$, as follows:
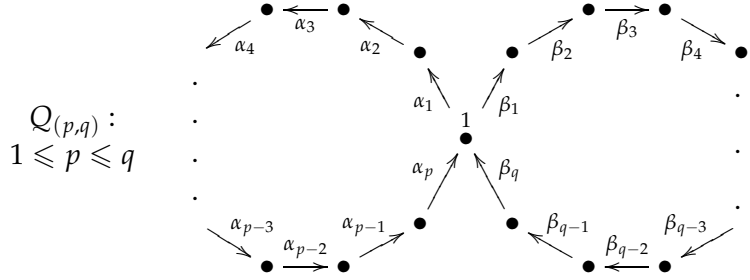
\[\xymatrix@C=.4pc@R=.4pc{
    &&&
     \bullet \ar^{\alpha_4}[lld] && \bullet \ar^{\alpha_3}[ll] &&& &&& \bullet \ar_{\beta_3}[rr] && \bullet
     \ar_{\beta_4}[rrd] \\
    &
  \ar@{*{}*{\rule[2.5mm]{0mm}{0mm}.\rule[-2.5mm]{0mm}{0mm}}*{}}[dddd]
      &&&&&&
      \bullet \ar^{\alpha_2}[ull] && \bullet \ar_{\beta_2}[urr] &&&&&& \bullet
  \ar@{*{}*{\rule[2.5mm]{0mm}{0mm}.\rule[-2.5mm]{0mm}{0mm}}*{}}[dddd]
      \\
    &&&&&&&&&&&&&&&&
    \\
     \save[] +<-3pc,0pt> *\txt<4pc>{$Q_{(p,q)}:$\\ $1\ppq p \ppq q$} \restore
      &&&&&&&& \stackrel{1}{\bullet} \ar^{\alpha_1}[luu] \ar_{\beta_1}[ruu]
      \\
     &&&&&&&&&&&&&&&&
     \\
    &
      \ar^{\alpha_{p-3}}[rrd] &&&&&&
        \bullet \ar^{\alpha_{p}}[ruu] && \bullet \ar_{\beta_{q}}[luu]
        &&&&&& \ar_{\beta_{q-3}}[lld] \\
    &&& \bullet \ar^{\alpha_{p-2}}[rr] && \bullet \ar^{\alpha_{p-1}}[rru]
      &&& &&&
      \bullet \ar_{\beta_{q-1}}[llu] && \bullet \ar_{\beta_{q-2}}[ll] \\
}\]

We denote the trivial path at the vertex $i$ by $e_i$. Paths are written from left to right. We write $\mo(\alpha)$ for
the trivial path corresponding to the origin of the arrow $\alpha$ and write $\mt(\alpha)$ for the trivial path
corresponding to the terminus of the arrow $\alpha$. The vertices of the quiver $\Q_{(p,q)}$ are labelled by $1, \ldots
, p+q-1$, in such a way that $\mo(\alpha_i) = i$ for $i = 1, \ldots , p$, and $\mt(\beta_j) = p+j$ for $j = 1, \ldots ,
q-1$. Thus $\mt(\alpha_i) = i+1$ for $i = 1, \ldots , p-1$, $\mt(\alpha_p) = 1$, $\mo(\beta_1) = 1$ and $\mo(\beta_j) =
p+j-1$ for $j = 2, \ldots , q$.

Set $\gamma=\alpha_1\alpha_{2}\cdots\alpha_p$ and $\delta=\beta_1\beta_{2}\cdots \beta_q.$ We define the following two admissible ideals in $K\Q_{(p,q)}$: \begin{enumerate}[\bf(a)] \item for a positive integer
$r$, let $I_r$ be the ideal generated by \[ \begin{array}{l} \alpha_p\alpha_1,\qquad \beta_q\beta_1,\qquad
(\gamma\delta)^r-(\delta\gamma)^r,\\ \alpha_i\cdots\alpha_p(\delta\gamma)^{r-1}\delta\alpha_1\cdots\alpha_{i} \text{
for all $2\ppq i\ppq p-1,$}\\ \beta_j\cdots\beta_q(\gamma\delta)^{r-1}\gamma\beta_1\cdots\beta_{j} \text{ for all $
2\ppq j\ppq q-1$;} \end{array}
 \]
\item for a pair of positive integers $(s,t)$, let $J_{(s,t)}$ be the ideal generated by \[ \begin{array}{l}
    \alpha_p\beta_1,\qquad \beta_q\alpha_1,\qquad \gamma^s-\delta^t, \\
    \alpha_i\cdots\alpha_p\gamma^{s-1}\alpha_1\cdots \alpha_i  \text{ for all $2\ppq i\ppq p-1$,}\\ \beta_j\cdots
    \beta_q \delta^{t-1}\beta_1\cdots\beta_j  \text{ for all $2\ppq j\ppq q-1,$} \end{array}
 \]
where, if $p=1$ then $s \pgq 2$, and, if $q=1$ then $p = 1, s\pgq 2$ and $t \pgq 2$.\\ \end{enumerate}

The algebras considered in \cite{BHS} are special cases of these algebras. Specifically $A(p,q) = K\Q_{(p,q)}/I_1$ so
that $r=1$, and $\Lambda(n) = K\Q_{(1,n)}/J_{(2,2)}$ so that $p=1, q=n$ and $s=2=t$. Moreover, some of these algebras are derived equivalent to algebras of dihedral type (see \cite{E}) in the classification of
Holm \cite{H}:  $K\Q_{(1,1)}/I_r=
D(1\A)^{r}_1$, $K\Q_{(1,2)}/I_r$ that is derived equivalent to $
D(2\B)^{1,r}(0)$, and $K\Q_{(2,2)}/I_r$ that is derived equivalent to $
D(3\K)^{r,1,1}$, all three of which come from tame blocks of finite groups  when $\car(K)=2$ and $r$ is a power of $2$, as well as $K\Q_{(2,2)}/J_{(s,t)}$ that is derived equivalent to $
D(2\R)^{1,s,t,1}$ and which does not come from blocks, see \cite{H1,H,L}.

\section{Background}

The following result and especially its consequences will be used repeatedly. They are given in \cite{AG}, but we
include the proofs here for completeness.

\bprop\label{prop:cycle}
Let $\mathcal Q$ be a quiver and let $I$ be an admissible ideal in $K{\mathcal Q}$ such that
$A=K{\mathcal Q}/I$ is a symmetric algebra. Let $\rho$ be a path in $\mathcal Q$ with $\rho\neq 0$ in $A$. Then there exists a
cycle $\rho\rho_1$ in $\mathcal Q$ with $\rho\rho_1$ and $\rho_1\rho$ non-zero in $A.$
\eprop

\begin{proof} Since $A$ is a symmetric algebra there exists a symmetric form $f:A\rightarrow K$ on $A$ whose kernel
contains no non-zero left or right ideals of $A.$ Then $\rho A$ is not contained in $\ker f$ so there exists a path
$\rho_1$ such that $f(\rho\rho_1)\neq 0.$ In particular, $\rho\rho_1\neq 0$ and $\mt(\rho)=\mo(\rho_1).$ Moreover,
since $f$ is symmetric, $f(\rho_1\rho)=f(\rho\rho_1)\neq 0$ so $\rho_1\rho\neq 0.$ Therefore $\mt(\rho_1)=\mo(\rho)$.
Hence $\rho\rho_1$ and $\rho_1\rho$ are cycles in $\mathcal Q$ which are non-zero in $A$. \end{proof}

\bc\label{cor:relations}
Let $A=K{\mathcal Q}/I$ be an indecomposable finite-dimensional symmetric special biserial algebra
that is not isomorphic to $N_1^1\cong K[X]/(X^2)$. Then:
\begin{enumerate}[(1)]
\item\label{uniquenonzerocomposition}
For any arrow $\alpha$, there is a unique arrow $\alpha'$ such that
    $\alpha'\alpha\neq 0$ and a unique arrow $\alpha''$ such that $\alpha\alpha''\neq 0.$
\item\label{degreevertices} For any vertex $v$ in $\mathcal Q$, the number of arrows that start at $v$ is equal to the
    number of arrows that end at $v$, and this number is either $1$ or $2.$
\end{enumerate} \ec

\begin{proof} The second statement follows easily from the first and the definition of a special biserial algebra. Here we prove the first statement.

First suppose that $\alpha$ is an arrow which is not a loop. Then $\alpha$ is a non-zero path, so, by Proposition
\ref{prop:cycle}, there exists a path $\rho\neq \alpha$ such that $\alpha\rho$ and $\rho\alpha$ are non-zero cycles.
Therefore we can take $\alpha'$ to be the last arrow in $\rho$ and $\alpha''$ the first arrow in $\rho.$ The uniqueness
of these arrows follows from the definition of a special biserial algebra.

Now suppose that $\alpha$ is a loop. Assume for contradiction that $\alpha\beta=0$ for every arrow $\beta$ in $\mathcal Q$.
Then $\alpha$ is in the socle of the indecomposable projective module $\mo(\alpha)A.$ If no other arrow starts at
$\mo(\alpha)$ then, since $A$ is indecomposable and $\alpha^2=0$, we get $A\cong K[X]/(X^2),$ a contradiction.
Therefore there is another arrow $\rho$ with $\mo(\alpha)=\mo(\rho).$ Choose a path $\sigma$ which is maximal with the
property $\rho\sigma\neq 0$  so that $\rho\sigma$ is in the socle of $\mo(\alpha)A.$ Since $A$ is a self-injective
algebra, $\soc(\mo(\alpha)A)$ is one-dimensional so that there exists a non-zero $c\in K$ such that
$\alpha=c\rho\sigma$, which is a contradiction since the ideal $I$ is admissible. Therefore there exists an arrow
$\alpha'$ with $\alpha'\alpha\neq0.$ The proof of the existence of the other arrow $\alpha''$ is similar. \end{proof}

\bigskip

\section{Classification theorem}\label{sec:quiver}

Our main theorem in this section is Theorem~\ref{thm:classification by quiver and relations} where we classify, by quiver and relations, all basic indecomposable finite-dimensional symmetric special biserial algebras with at most one
non-uniserial indecomposable projective module.

\bprop \label{prop:algebrassatisfyproperty}
The algebras $K\Q_{(p,q)}/I_r$ and $K\Q_{(p,q)}/J_{(s,t)}$ are symmetric special
biserial algebras with at most one non-uniserial indecomposable projective module.
\eprop

\begin{proof} It is easy to see that these algebras are special biserial and that all but one of each of their
indecomposable projective modules are uniserial. Moreover, the algebras are weakly symmetric, that is, the top and the
socle of each indecomposable projective module are isomorphic. It remains to prove that the algebras are symmetric.

For $A=K\Q_{(p,q)}/I_r$, the socle of $A$ has a $K$-basis consisting of $(\delta\gamma)^r,$
$\alpha_i\cdots\alpha_p(\delta\gamma)^{r-1}\alpha_1\cdots\alpha_{i-1}$ and
$\beta_j\cdots\beta_q(\gamma\delta)^{r-1}\beta_1\cdots\beta_{j-1}$ for $2\ppq i\ppq p$ and $2\ppq j\ppq q,$ that is,
all the paths obtained from cyclic permutations of $(\gamma\delta)^r,$ where we recall that
$\gamma=\alpha_1\alpha_{2}\cdots\alpha_p$ and $\delta=\beta_1\beta_{2}\cdots \beta_q$. Complete this $K$-basis of
$\soc(A)$ with paths in $K\Q_{(p,q)}$ to obtain a basis of $A$, and define $f:A\rightarrow K$ on this basis by sending
the elements in the socle to $1$ and the others to $0.$ Then it follows from \cite[Proposition 3.1]{HZ}, that $\ker f$
contains no non-zero left or right ideals of $A$. Moreover, $f$ is clearly  symmetric, since the paths on which it is
non-zero are all the cyclic permutations of a single path. Thus $A$ is a symmetric algebra.

For $A=K\Q_{(p,q)}/I_{(s,t)}$ the argument is similar, but this time the socle of $A$ is generated by all the cyclic
permutations of the two paths $\gamma^s$ and $\delta^t.$ \end{proof}

We now have the following result.

\bt\label{thm:classification by quiver and relations}
Let $A$ be a basic indecomposable finite-dimensional symmetric special biserial algebra with at most one
non-uniserial indecomposable projective module. Then $A$ is isomorphic to a Nakayama algebra $N_m^n$ for some positive
integers $m$ and $n$, or to $K\Q_{(p,q)}/I_r$ for some positive integers $p,$ $q$ and $r$, or to
$K\Q_{(p,q)}/J_{(s,t)}$ for some positive integers $p,$ $q$, $s$ and $t.$
\et

\begin{proof} Set $A=K{\mathcal Q}/I$ for some quiver $\mathcal Q$ and some admissible ideal $I.$  It is already known that a
symmetric (special biserial) algebra with no non-uniserial indecomposable projective module is isomorphic to a Nakayama
algebra $N_m^n$. We may therefore assume that all except one indecomposable projective $A$-modules are uniserial.
Consequently, using Corollary \ref{cor:relations}(2), we must have one vertex that is the end point of exactly two
arrows and the starting point of exactly two arrows, which we label $1$, and the other vertices are the end point of
exactly one arrow and the starting point of exactly one arrow. Therefore the quiver of $A$ must be $\Q_{(p,q)}$ for
some positive integers $p,q.$ Without loss of generality we may assume that $p\ppq q.$

Now consider the composition $\alpha_p\alpha_1.$ There are two cases: $\alpha_p\alpha_1=0$ and $\alpha_p\alpha_1\neq
0.$

First assume that $\alpha_p\alpha_1=0.$ Then it follows from Corollary \ref{cor:relations}(1) that $\alpha_p\beta_1\neq
0,$ $\beta_q\alpha_1\neq0$ and $\beta_q\beta_1=0.$ Now, for each vertex $k$ with $k \neq 1$, $e_kA$ is uniserial and,
since $A$ is symmetric, $\top{(e_kA)}\cong \soc{(e_kA)}\cong S_k$, the simple module at $k$. Therefore, for $i\neq 1,$
we get a relation $\alpha_i\cdots\alpha_p(\delta\gamma)^{u_i}\delta\alpha_1\cdots\alpha_i=0$ with
$\alpha_i\cdots\alpha_p(\delta\gamma)^{u_i}\delta\alpha_1\cdots\alpha_{i-1}\neq0$ for some integer $u_i$, and, for
$j\neq 1,$ we get a relation $\beta_j\cdots\beta_p(\gamma\delta)^{v_j}\gamma\beta_1\cdots\beta_j=0$ with
$\beta_j\cdots\beta_p(\gamma\delta)^{v_j}\gamma\beta_1\cdots\beta_{j-1}\neq0$ for some integer $v_j.$ Now consider
$e_1A.$ Since $\rad(e_1A) = \alpha_1A + \beta_1A$ and $\soc(e_1A) \cong S_1$, there is an element in $\soc(e_1A)$ of
the form $(\gamma\delta)^r$ or $(\gamma\delta)^r\gamma$ for some integer $r$, and there is an element of $\soc(e_1A)$
of the form $(\delta\gamma)^s$ or $(\delta\gamma)^s\delta$ for some integer $s$. But, $\soc(e_1A)$ is simple so we must have a relation of one of the following forms:
\begin{enumerate}[(i)]
\item \sloppy $(\gamma\delta)^r\gamma=c(\delta\gamma)^s\neq 0$ for some non-zero $c\in K.$ Note that since $I$ is
    admissible we must have $s\pgq 1.$ But then, if $r\pgq s,$ we would have
    $(\gamma\delta)^r\gamma=\gamma(\delta\gamma)^s(\delta\gamma)^{r-s} = c^{-1}\gamma(\gamma\delta)^r\gamma(\delta\gamma)^{r-s}=0$,
    which is a contradiction, and if $r<s$ we would have
    $(\delta\gamma)^s=\delta(\gamma\delta)^{s-1-r}(\gamma\delta)^r\gamma=c\delta(\gamma\delta)^{s-1-r}(\delta\gamma)^s = 0$,
    which is also a contradiction. Therefore we cannot have this type of relation.
\item $(\gamma\delta)^r=c(\delta\gamma)^s\delta$ for some non-zero $c\in K.$ As in the previous case, this relation
    cannot occur.
\item $(\gamma\delta)^r\gamma=c(\delta\gamma)^s\delta\neq 0$  for some non-zero $c\in K.$ Then, if $r>s,$ we have
    $(\gamma\delta)^r\gamma=\gamma(\delta\gamma)^r=\gamma(\delta\gamma)^s\delta(\gamma\delta)^{r-1-s}\gamma = c^{-1}\gamma(\gamma\delta)^r\gamma(\gamma\delta)^{r-1-s}\gamma=0,$
    which is a contradiction, and if $s>r$ we have a similar contradiction. Therefore $s=r.$ Now we also have
    $\alpha_2\cdots\alpha_p(\delta\gamma)^{u_2}\delta\alpha_1\alpha_2=0$ so, multiplying on the left by $\alpha_1$
    and on the right by $\alpha_3\cdots \alpha_p$, we get $\gamma(\delta\gamma)^{u_2+1}=0$ so that $u_2+1>r.$ But
    $\alpha_2\cdots
    \alpha_p(\delta\gamma)^r\delta\alpha_1=c^{-1}\alpha_2\cdots\alpha_p(\gamma\delta)^r\gamma\alpha_1=0$ so $u_2<r$
    which is impossible. Therefore a relation of this form cannot occur either.
\item $(\gamma\delta)^r=c(\delta\gamma)^s$ for some non-zero $c\in K$, and this is the only possible type of
    relation. Here again, if $r>s,$ then
    $(\gamma\delta)^r=\gamma(\delta\gamma)^s(\delta\gamma)^{r-s-1}\delta = c^{-1}\gamma(\gamma\delta)^r(\delta\gamma)^{r-s-1}\delta=0$
    which is a contradiction and if $s>r$ we get a similar contradiction. Therefore $r=s$ so that the relation is
    $(\gamma\delta)^r=c(\delta\gamma)^r$ for some $r\pgq 1$ and $c\in K^*$.
\end{enumerate}

\sloppy Given this relation, we are now able to determine the $u_i$ and the $v_j.$ Since $\alpha_i\cdots
\alpha_p(\delta\gamma)^r\delta\alpha_1\cdots\alpha_{i-1}$ $=
c^{-1}\alpha_i\cdots\alpha_p(\gamma\delta)^r\delta\alpha_1\cdots\alpha_{i-1}=0,$ we must have $r>u_i.$ Moreover,
$$(\gamma\delta)^{u_i+2}=\alpha_1\cdots
\alpha_{i-1}(\alpha_i\cdots\alpha_p(\delta\gamma)^{u_i}\delta\alpha_1\cdots\alpha_i)\alpha_{i+1}\cdots\alpha_p\delta=0$$
so that $u_i+2>r.$ Hence $u_i=r-1$ for all $i = 2, \ldots , p.$ Similarly, $v_j=r-1$ for all $j = 2, \ldots , q$.
Moreover, we note that the relations $\alpha_p(\delta\gamma)^{r-1}\delta\alpha_1\cdots\alpha_p=0$ (when $i=p$) and
$\beta_p(\gamma\delta)^{r-1}\gamma\beta_1\cdots\beta_p=0$ (when $j = q$) are superfluous, so are not required to give a minimal generating set of the ideal $I_r$.

Finally, we show that $c$ must be equal to $1.$ Since $A$ is symmetric, there exists a symmetric linear map
$f:A\rightarrow K$ whose kernel does not contain any non-zero left or right ideal. In particular, the socle of $A$ is
not contained in $\ker f.$ But, from the relations obtained above, we see that the socle is generated as a $K$-vector
space by all the paths obtained by cyclic permutations of $(\gamma\delta)^r$. Since $f$ is symmetric it follows that
$f((\gamma\delta)^r)\neq 0.$ But we have \[
f((\gamma\delta)^r)=f((\delta\gamma)^r)=f(c(\gamma\delta)^r)=cf((\gamma\delta)^r) \] so that $c=1.$

Hence we have shown that  $A\cong K\Q_{(p,q)}/I_r.$

\smallskip

Now assume that $\alpha_p\alpha_1\neq 0.$ Then it follows from Corollary \ref{cor:relations}(1) that $\alpha_p\beta_1=
0,$ $\beta_q\alpha_1=0$ and $\beta_q\beta_1\neq0.$ The same methods as in the previous case show that we must have a
relation of the form $\gamma^s=c\delta^t$ for some non-zero $c\in K$ and some positive integers $s$ and $t$ (by
considering the structure of $e_1A$) and relations $\alpha_i\cdots\alpha_p\gamma^{s-1}\alpha_1\cdots \alpha_i$ for
all $2\ppq i\ppq p$ and $\beta_j\cdots \beta_q \delta^{t-1}\beta_1\cdots\beta_j$ for all $2\ppq j\ppq q$ (from the
structure of the other indecomposable projectives and using the relation $\gamma^s=c\delta^t$). Moreover, since $K$ is
algebraically closed, we may replace $\alpha_1$ by $c'\alpha_1$, where $c'$ is a $t$-th root of $c$, and thus we may
replace the relation $\gamma^s=c\delta^t$ by $\gamma^s=\delta^t$. Again we may find a minimal set of relations, and so
conclude that $A\cong K\Q_{(p,q)}/J_{(s,t)}.$ \end{proof}

\bigskip

\section{Hochschild cohomology groups}\label{sec:HH}

Our aim is now to investigate the derived equivalences among these algebras. It is well-known that Hochschild
cohomology is an invariant under derived equivalence, and this section determines some of the Hochschild cohomology
groups of the algebras $K{\mathcal Q}_{(p,q)}/I_r$ and $K{\mathcal Q}_{(p,q)}/J_{(s,t)}$, so that a full classification up to derived equivalence can be given in Section \ref{sec:final}.

\sloppy Let $\Gamma(p,q;r) = K{\mathcal Q}_{(p,q)}/I_r$ and $\Lambda(p,q;s,t) = K{\mathcal Q}_{(p,q)}/J_{(s,t)}$.
The special cases of the algebras $\Gamma(p,q;1)$ and $\Lambda(1,n;2,2)$ were considered in \cite{BHS}, where, in their notation, we have $A(p,q) = K{\mathcal Q}_{(p,q)}/I_1 = \Gamma(p,q;1)$ and $\Lambda(n) = K{\mathcal Q}_{(1,n)}/ J_{(2,2)} = \Lambda(1,n;2,2)$.

We begin by describing $\HH^0(A)$ and $\HH^1(A)$ for the algebras $A = \Lambda(p,q;s,t)$ and $A = \Gamma(p,q;r)$. Recall that $\HH^0(A) = Z(A)$, the centre of the algebra $A$.

\bigskip

\subsection{$\HH^0(A)$ and $\HH^1(A)$ for the algebra $A=\Lambda(p,q;s,t)$}

We begin with the algebra $\Lambda(p,q;s,t) = K{\mathcal Q}_{(p,q)}/J_{(s,t)}$ where $1 \ppq p \ppq q$. Recall that
$\gamma = \alpha_1\alpha_2\cdots\alpha_p$ and $\delta = \beta_1\beta_2\cdots\beta_q$. Let $\gamma_i =
\alpha_i\cdots\alpha_p\alpha_1\cdots\alpha_{i-1}$ for $i = 1, \ldots , p$, and let $\delta_j =
\beta_j\cdots\beta_q\beta_1\cdots\beta_{j-1}$ for $j = 1, \ldots , q$, so that $\gamma = \gamma_1$ and $\delta =
\delta_1$.

\bprop\label{prop:centre of lambda}
Consider the algebra $\Lambda(p,q;s,t)$ and let $1 \ppq p \ppq q$. Let $x =
\sum_{i=1}^p \gamma_i$ and $y = \sum_{j=1}^q \delta_j$. Then $$\dim\HH^0(\Lambda(p,q;s,t)) = p+q+s+t-2$$ and the set
$$\set{1, x, \ldots , x^{s-1}, y, \ldots , y^{t-1}, \gamma_i^s, \delta_j^t \mbox{ for } i = 1, \ldots, p \mbox{ and } j
= 2, \ldots , q}$$ is a $K$-basis of $\HH^0(\Lambda(p,q;s,t))$.
\eprop

\begin{proof} We note that $\gamma_i\beta_j = 0 = \beta_j\gamma_i$ and $\delta_j\alpha_i = 0 = \alpha_i\delta_j$ for
all $i = 1, \ldots, p, \ j = 1, \ldots , q$ so $x, y \in Z(\Lambda(p,q;s,t))$. Moreover $\Lambda(p,q;s,t)$ is weakly
symmetric so all socle elements are central, namely, $\gamma_i^s$ and $\delta_j^t$ are central for $i = 1, \ldots, p$,
$j = 2, \ldots , q$. The result now follows. \end{proof}

We remark that Proposition \ref{prop:centre of lambda} may be simplified if $p = q = 1$, as then $\Lambda(1,1;s,t)$ is
the commutative algebra $K[\alpha,\beta]/(\alpha\beta, \alpha^s-\beta^t)$ with $s \pgq 2, t \pgq 2$. Thus
$\HH^0(\Lambda(1,1;s,t)) = \Lambda(1,1;s,t)$ which has $K$-basis $\set{1, \alpha, \ldots, \alpha^s, \beta, \ldots ,
\beta^{t-1}}$ and dimension $s+t$.

\bigskip

In order to compute the first Hochschild cohomology group, we now fix a minimal set of generators of the ideal
$J_{(s,t)}$, and denote this set by $g^2$.

\bprop\label{prop:min gens for lambda}
Consider the algebra $\Lambda(p,q;s,t)$ with $1 \ppq p \ppq q$.

If $p \pgq 2$ then the following elements form a minimal set of generators for the ideal $J_{(s,t)}$:
$$\begin{array}{lcl} g^2_1 & = & \gamma^s - \delta^t\\ g^2_i & = & \alpha_i\cdots\alpha_p\gamma^{s-1}\alpha_1\cdots
\alpha_i  \text{ for all $2\ppq i\ppq p-1$}\\ g^2_p & = & \alpha_p\beta_1\\ g^2_{p+1} & = & \beta_q\alpha_1\\ g^2_{p+j}
& = & \beta_j\cdots \beta_q \delta^{t-1}\beta_1\cdots\beta_j  \text{ for all $2\ppq j\ppq q-1$}.\\ \end{array}$$

If $p = 1$ then the following elements form a minimal set of generators for the ideal $J_{(s,t)}$: $$\begin{array}{lcl}
g^2_0 & = & \alpha^s - \delta^t\\ g^2_1 & = & \alpha\beta_1\\ g^2_2 & = & \beta_q\alpha\\ g^2_{j+1} & = & \beta_j\cdots
\beta_q \delta^{t-1}\beta_1\cdots\beta_j  \text{ for all $2\ppq j\ppq q-1$.}\\ \end{array}$$
\eprop

\bigskip

We now compute the first Hochschild cohomology group $\HH^1(\Lambda(p,q;s,t))$. We use the explicit description of the
start of a minimal projective bimodule resolution $(P^*, d^*)$ for $\Lambda(p,q;s,t)$ as given in \cite{GSn}. All
tensors are over $K$ so we write $\otimes$ for $\otimes_K$. For ease of notation, write $\Lambda$ for $\Lambda(p,q;s,t)$.
Let $P^2 = \bigoplus_{k}\Lambda \mo(g^2_k)\otimes \mt(g^2_k)\Lambda, \ P^1 = \bigoplus_{a\ {\rm arrow}}\Lambda
\mo(a)\otimes \mt(a)\Lambda$ and $P^0 = \bigoplus_{v=1}^{p+q-1}\Lambda e_v \otimes e_v\Lambda$. Then the minimal
projective bimodule resolution of $\Lambda$ begins $$\cdots \longrightarrow P^2 \stackrel{d^2}{\longrightarrow} P^1
\stackrel{d^1}{\longrightarrow} P^0  \stackrel{d^0}{\longrightarrow} \Lambda \longrightarrow 0$$ with the following
maps. The map $d^0$ is the usual multiplication map. The map $d^1 \colon P^1 \to P^0$ is given by $$d^1 \colon
\mo(a)\otimes \mt(a) \mapsto \mo(a)\otimes a - a\otimes \mt(a)$$ where the first term $\mo(a)\otimes a$ lies in the
summand $\Lambda \mo(a) \otimes \mo(a)\Lambda$ and the second term $a\otimes \mt(a)$ lies in the summand $\Lambda \mt(a)
\otimes \mt(a)\Lambda$. Now, each element of $g^2$ is a linear combination of paths in $K{\mathcal Q}_{(p,q)}$ so, for $x \in
g^2$, we may write $$x = \sum_{j=1}^{r}c_j a_{1 j} \ldots a_{k j} \ldots a_{s_j j}$$ where $c_j \in K$ and the $a_{k
j}$ are arrows in ${\mathcal Q}_{(p,q)}$. With this notation for $x \in g^2$, the map $d^2 \colon P^2 \to P^1$ is given by
$$d^2 \colon \mo(x) \otimes \mt(x) \mapsto \sum_{j=1}^{r}c_j \sum_{k=1}^{s_j}a_{1 j} \ldots a_{k-1 j} \otimes a_{k+1 j}
\ldots a_{s_j j}$$ where the term $a_{1 j} \ldots a_{k-1 j} \otimes a_{k+1 j} \ldots a_{s_j j}$ lies in the summand
$\Lambda\mo(a_{k j}) \otimes \mt(a_{k j})\Lambda$ of $P^1$.

\bigskip

We now apply $\Hom_{\Lambda^e}(-, -)$ to this resolution, where $\Lambda^e=\Lambda\ot \Lambda^{op}$ is the enveloping
algebra of $\Lambda$. Let $\partial^1 : \Hom_{\Lambda^e}(P^1, \Lambda) \to \Hom_{\Lambda^e}(P^2, \Lambda)$ be the map
induced by $d^2$ and let $\partial^0 : \Hom_{\Lambda^e}(P^0, \Lambda) \to \Hom_{\Lambda^e}(P^1, \Lambda)$ be the map
induced by $d^1$. Then $\HH^1(\Lambda) = \ker\partial^1/\im\partial^0$.

\bprop\label{prop:HH1 for lambda} If $q \pgq 2$ then $$\dim\HH^1(\Lambda(p,q;s,t)) = \left \{ \begin{array}{ll} s+t &
\mbox{ if $\car K \mid \gcd(s,t)$}\\ s+t-1 & \mbox{ otherwise}. \end{array} \right. $$ If $q = 1$ then $p=1$ and
$$\dim\HH^1(\Lambda(1,1;s,t)) = \left \{ \begin{array}{ll} s+t+1 & \mbox{ if $\car K \mid \gcd(s,t)$}\\ s+t & \mbox{
otherwise}. \end{array} \right. $$ \eprop

\begin{proof} There are three cases to consider.

\noindent {\bf Case I}: $p \pgq 2$

We start by calculating $\im\partial^0$. Let $\varphi \in \Hom(P^0, \Lambda(p,q;s,t))$ so that $\partial^0(\varphi) =
\varphi d^1$. Suppose that $\varphi$ is given by $$\begin{array}{ccl} \varphi: & e_1\otimes e_1 & \mapsto \ c_{1,0}e_1
+ c_{1,1}\gamma_1 + \cdots + c_{1,s}\gamma_1^s + d_{1,1}\delta_1 + \cdots + d_{1,t-1}\delta_1^{t-1}\\ & e_i \otimes e_i
& \mapsto \ c_{i,0}e_i + c_{i,1}\gamma_i + \cdots + c_{i,s}\gamma_i^{s} \mbox{\ \ for $i = 2, \ldots , p$}\\ &
e_{p-1+i} \otimes e_{p-1+i} & \mapsto \ d_{i,0}e_{p-1+i} + d_{i,1}\delta_i + \cdots + d_{i,t}\delta_i^{t} \mbox{\ \ for
$i = 2, \ldots , q$} \end{array}$$ where $c_{i,j}, d_{i,j} \in K$.

We have $$\begin{array}{ccl} \varphi d^1(\mo(\alpha_1)\otimes \mt(\alpha_1)) & = & \varphi(e_1\otimes \alpha_1 -
\alpha_1\otimes e_2) \ = \ \varphi(e_1\otimes e_1)\alpha_1 - \alpha_1\varphi(e_2\otimes e_2)\\ & = & (c_{1,0} -
c_{2,0})\alpha_1 + (c_{1,1} - c_{2,1})\gamma_1\alpha_1 + \cdots + (c_{1,s-1} - c_{2,s-1})\gamma_1^{s-1}\alpha_1.
\end{array}$$ In a similar way we get $$\begin{array}{ccl} \varphi d^1(\mo(\alpha_2)\otimes \mt(\alpha_2)) & = &
\varphi(e_2\otimes \alpha_2 - \alpha_2\otimes e_3)\\ & = & (c_{2,0}-c_{3,0})\alpha_2 + (c_{2,1} -
c_{3,1})\gamma_2\alpha_2 + \cdots + (c_{2,s-1} - c_{3,s-1})\gamma_2^{s-1}\alpha_2\\ \\ \vdots &  &  \\ \\ \varphi
d^1(\mo(\alpha_p)\otimes \mt(\alpha_p)) & = & \varphi(e_p\otimes \alpha_p - \alpha_p\otimes e_1)\\ & = &
(c_{p,0}-c_{1,0})\alpha_p + (c_{p,1} - c_{1,1})\gamma_p\alpha_p + \cdots + (c_{p,s-1} -
c_{1,s-1})\gamma_p^{s-1}\alpha_p\\ \\ \varphi d^1(\mo(\beta_1)\otimes \mt(\beta_1)) & = & \varphi(e_1\otimes \beta_1 -
\beta_1\otimes e_{p+1})\\ & = & (c_{1,0} - d_{2,0})\beta_1 + (d_{1,1} - d_{2,1})\delta_1\beta_1 + \cdots + (d_{1,t-1} -
d_{2,t-1})\delta_1^{t-1}\beta_1\\ \\ \vdots &  &  \\ \\ \varphi d^1(\mo(\beta_q)\otimes \mt(\beta_q)) & = &
\varphi(e_{p+q-1}\otimes \beta_q - \beta_q\otimes e_1)\\ & = & (d_{q,0}-c_{1,0})\beta_q + (d_{q,1} -
d_{1,1})\delta_q\beta_q + \cdots + (d_{q,t-1} - d_{1,t-1})\delta_q^{t-1}\beta_q. \end{array}$$ Thus $\dim\im\partial^0
= (p-1)s+(q-1)t$.

\medskip

Now let $\psi \in \ker\partial^1$ so that $\psi d^2 = 0$, and suppose that $\psi \in \Hom(P^1, \Lambda(p,q;s,t))$ is
given by $$\begin{array}{ccl} \psi: & \mo(\alpha_i)\otimes \mt(\alpha_i) & \mapsto \ c_{i,0}\alpha_i +
c_{i,1}\gamma_i\alpha_i + \cdots + c_{i,s-1}\gamma_i^{s-1}\alpha_i\\ & \mo(\beta_j)\otimes \mt(\beta_j) & \mapsto \
d_{j,0}\beta_j + d_{j,1}\delta_j\beta_j + \cdots + d_{j,t-1}\delta_j^{t-1}\beta_j \end{array}$$ for $i = 1, \dots, p,\
j = 1, \dots, q$ and where $c_{k,l}, d_{k,l} \in K$.

From Proposition \ref{prop:min gens for lambda} with $p \pgq 2$, and recalling that $\gamma = \gamma_1$ and $\delta = \delta_1$, it is easy to see that $\psi d^2(\mo(g^2_k)\otimes \mt(g^2_k))$ is immediately zero for $k = 2, \dots , p+q-1$ and so
we do not get any restrictions on the constants in the cases where $g^2_k$ is a monomial. It remains to consider $g^2_1
= \gamma^s - \delta^t$. The condition $\psi d^2(\mo(g^2_1)\otimes \mt(g^2_1)) = 0$ gives that $$(s(c_{1,0} + c_{2,0} +
\cdots + c_{p,0}) - t(d_{1,0} + d_{2,0} + \cdots + d_{q,0}))\gamma^s = 0$$ so that $s(c_{1,0} + c_{2,0} + \cdots +
c_{p,0}) - t(d_{1,0} + d_{2,0} + \cdots + d_{q,0}) = 0$. Hence $$\dim\ker\partial^1 = \left \{ \begin{array}{ll} ps+qt
& \mbox{ if $\car K \mid \gcd(s,t)$}\\ ps+qt-1 & \mbox{ otherwise}. \end{array} \right. $$

Thus, for $p \pgq 2$, we have $$\dim\HH^1(\Lambda(p,q;s,t)) = \left \{ \begin{array}{ll} s+t & \mbox{ if $\car K \mid
\gcd(s,t)$}\\ s+t-1 & \mbox{ otherwise}. \end{array} \right. $$

\bigskip

\noindent {\bf Case II}: $p = 1$ and $q \pgq 2$

To calculate $\im\partial^0$, let $\varphi \in \Hom(P^0, \Lambda(1,q;s,t))$ so that $\partial^0(\varphi) = \varphi d^1$.
Suppose that \linebreak $\varphi \in \Hom(P^0, \Lambda(1,q;s,t))$ is given by $$\begin{array}{ccl} \varphi: & e_1\otimes
e_1 & \mapsto \ c_{1,0}e_1 + c_{1,1}\alpha + \cdots + c_{1,s}\alpha^s + d_{1,1}\delta_1 + \cdots +
d_{1,t-1}\delta_1^{t-1}\\ & e_i \otimes e_i & \mapsto \ d_{i,0}e_i + d_{i,1}\delta_i + \cdots + d_{i,t}\delta_i^{t}
\end{array}$$ for $i = 2, \dots, q$ and where $c_{1,j}, d_{1,j}, d_{i,j} \in K$.

We have $$\varphi d^1(\mo(\alpha)\otimes \mt(\alpha)) = \varphi(e_1\otimes \alpha - \alpha\otimes e_1) =
\varphi(e_1\otimes e_1)\alpha - \alpha\varphi(e_1\otimes e_1) = 0$$ since $\alpha \in Z(\Lambda(1,q;s,t))$. Now,
$$\begin{array}{ccl} \varphi d^1(\mo(\beta_1)\otimes \mt(\beta_1)) & = & \varphi(e_1\otimes \beta_1 - \beta_1\otimes
e_2)\\ & = & (c_{1,0}-d_{2,0})\beta_1 + (d_{1,1}-d_{2,1})\delta_1\beta_1 + \cdots +
(d_{1,t-1}-d_{2,t-1})\delta_1^{t-1}\beta_1\\ \\ \vdots &  &  \\ \\ \varphi d^1(\mo(\beta_q)\otimes \mt(\beta_q)) & = &
\varphi(e_q\otimes \beta_q - \beta_q\otimes e_1)\\ & = & (d_{q,0} - c_{1,0})\beta_q + (d_{q,1} -
d_{1,1})\delta_q\beta_q + \cdots + (d_{q,t-1} - d_{1,t-1})\delta_q^{t-1}\beta_q.\\ \end{array}$$ Thus
$\dim\im\partial^0 = (q-1)t$.

\medskip

Now let $\psi \in \ker\partial^1$ so that $\psi d^2 = 0$, and suppose that $\psi \in \Hom(P^1, \Lambda(1,q;s,t))$ is
given by $$\begin{array}{ccl} \psi: & \mo(\alpha)\otimes \mt(\alpha) & \mapsto \ c_{1,0}e_1 + c_{1,1}\alpha + \cdots +
c_{1,s}\alpha^s + \tilde{c}_{1,1}\delta_1 + \cdots + \tilde{c}_{1,t-1}\delta_1^{t-1}\\ & \mo(\beta_i)\otimes
\mt(\beta_i) & \mapsto \ d_{i,0}\beta_i + d_{i,1}\delta_i\beta_i + \cdots + d_{i,t-1}\delta_i^{t-1}\beta_i
\end{array}$$ for $i = 1, \dots, q$ and where $c_{1,j}, \tilde{c}_{1,j}, d_{i,j} \in K$.

From Proposition \ref{prop:min gens for lambda}, the minimal generating set for $J_{(s,t)}$ is $$\{g^2_0 = \alpha^s -
\delta^t, \ g^2_1 = \alpha\beta_1, \ g^2_2 = \beta_q\alpha, \ g^2_{j+1} = \beta_j\cdots \beta_q
\delta^{t-1}\beta_1\cdots\beta_j \mid 2\ppq j\ppq q-1\}$$ where we recall that $\delta = \delta_1$.
Starting with $g^2_1$, the equation $\psi d^2(\mo(g^2_1)\otimes \mt(g^2_1)) = 0$ gives that $$\begin{array}{lcl} 0 & =
& \psi(\mo(\alpha)\otimes\mt(\alpha))\beta_1 + \alpha\psi(\mo(\beta_1)\otimes \mt(\beta_1))\\ & = & (c_{1,0}e_1 +
\tilde{c}_{1,1}\delta_1 + \cdots + \tilde{c}_{1,t-1}\delta_1^{t-1})\beta_1.\\ \end{array}$$ Hence $c_{1,0} =
\tilde{c}_{1,1} = \cdots = \tilde{c}_{1,t-1} = 0$. So we may immediately simplify our expression for $\psi$ as
$$\begin{array}{ccl} \psi: & \mo(\alpha)\otimes \mt(\alpha) & \mapsto \ c_{1,1}\alpha + \cdots + c_{1,s}\alpha^s \\ &
\mo(\beta_i)\otimes \mt(\beta_i) & \mapsto \ d_{i,0}\beta_i + d_{i,1}\delta_i\beta_i + \cdots +
d_{i,t-1}\delta_i^{t-1}\beta_i \end{array}$$ for $i = 1, \dots, q$. It then follows that $\psi d^2(\mo(g^2_j)\otimes
\mt(g^2_j))$ is zero for $j = 2, \dots , q$ and so we do not get any restrictions on the constants here.
Finally, $$0 = \psi d^2(\mo(g^2_0)\otimes \mt(g^2_0)) = (sc_{1,1} - t(d_{1,0} + d_{2,0} + \cdots + d_{q,0}))\alpha^s$$
so that $sc_{1,1} - t(d_{1,0} + d_{2,0} + \cdots + d_{q,0}) = 0$. Hence $$\dim\ker\partial^1 = \left \{
\begin{array}{ll} s+qt & \mbox{ if $\car K \mid \gcd(s,t)$}\\ s+qt-1 & \mbox{ otherwise}. \end{array} \right. $$

Thus, for $q \pgq 2$, we have $$\dim\HH^1(\Lambda(1,q;s,t)) = \left \{ \begin{array}{ll} s+t & \mbox{ if $\car K \mid
\gcd(s,t)$}\\ s+t-1 & \mbox{ otherwise}. \end{array} \right. $$

\bigskip

\noindent {\bf Case III}: $p = 1 = q$

To calculate $\im\partial^0$, let $\varphi \in \Hom(P^0, \Lambda(1,1;s,t))$ so that $\partial^0(\varphi) = \varphi d^1$.
We have $$\varphi d^1(\mo(\alpha)\otimes \mt(\alpha)) = \varphi(e_1\otimes \alpha - \alpha\otimes e_1) =
\varphi(e_1\otimes e_1)\alpha - \alpha\varphi(e_1\otimes e_1) = 0$$ since $\Lambda(1,1;s,t)$ is commutative. Similarly
$\varphi d^1(\mo(\beta)\otimes \mt(\beta)) = 0$. Thus $\im\partial^0 = (0)$.

Hence $\HH^1(\Lambda(1,1;s,t)) = \ker\partial^1$. Let $\psi \in \ker\partial^1$ so that $\psi d^2 = 0$, and suppose that
$\psi \in \Hom(P^1, \Lambda(1,1;s,t))$ is given by $$\begin{array}{ccl} \psi: & \mo(\alpha)\otimes \mt(\alpha) & \mapsto
\ c_{1,0}e_1 + c_{1,1}\alpha + \cdots + c_{1,s}\alpha^s + d_{1,1}\beta + \cdots + d_{1,t-1}\beta^{t-1}\\ &
\mo(\beta)\otimes \mt(\beta) & \mapsto \ c_{2,0}e_1 + c_{2,1}\alpha + \cdots + c_{2,s}\alpha^s + d_{2,1}\beta + \cdots
+ d_{2,t-1}\beta^{t-1} \end{array}$$ where $c_{i,j}, d_{i,j} \in K$.

From Proposition \ref{prop:min gens for lambda}, the minimal generating set for $J_{(s,t)}$ is $\{g^2_0 = \alpha^s -
\beta^t, \ g^2_1 = \alpha\beta, \ g^2_2 = \beta\alpha\}$. Starting with $g^2_1$, the equation $\psi
d^2(\mo(g^2_1)\otimes \mt(g^2_1)) = 0$ gives that $$c_{1,0}\beta + d_{1,1}\beta^2 + \cdots + d_{1,t-2}\beta^{t-1} +
c_{2,0}\alpha + c_{2,1}\alpha^2 + \cdots + c_{2,s-2}\alpha^{s-1} + (d_{1,t-1} +c_{2,s-1})\alpha^{s} = 0.$$ Hence
$c_{1,0} = d_{1,1} = \cdots = d_{1,t-2} = c_{2,0} = c_{2,1} = \cdots = c_{2,s-2} = 0$ and $d_{1,t-1} + c_{2,s-1} = 0$.
So we may simplify our expression for $\psi$ as $$\begin{array}{ccl} \psi : \mo(\alpha)\otimes \mt(\alpha)
& \mapsto & c_{1,1}\alpha + \cdots + c_{1,s}\alpha^s + d_{1,t-1}\beta^{t-1}\\ \psi : \mo(\beta)\otimes \mt(\beta) &
\mapsto & -d_{1,t-1}\alpha^{s-1} + c_{2,s}\alpha^s + d_{2,1}\beta + \cdots + d_{2,t-1}\beta^{t-1}. \end{array}$$ We
then have that $$0 = \psi d^2(\mo(g^2_0)\otimes \mt(g^2_0)) = sc_{1,1}\alpha^s - td_{2,1}\beta^t = (sc_{1,1} -
td_{2,1})\alpha^s$$ and hence $sc_{1,1} - td_{2,1} = 0$. The final equation $\psi d^2(\mo(g^2_0)\otimes \mt(g^2_0)) =
0$ gives no new information. Hence $\psi \in \ker\partial^1$ is given by $$\begin{array}{ccl} \psi : \mo(\alpha)\otimes
\mt(\alpha) & \mapsto & c_{1,1}\alpha + \cdots + c_{1,s}\alpha^s + d_{1,t-1}\beta^{t-1}\\ \psi : \mo(\beta)\otimes
\mt(\beta) & \mapsto &  -d_{1,t-1}\alpha^{s-1} + c_{2,s}\alpha^s + d_{2,1}\beta + \cdots + d_{2,t-1}\beta^{t-1}
\end{array}$$ with the additional linear dependency that $sc_{1,1} - td_{2,1} = 0$. Therefore $$\dim\HH^1(\Lambda(1,1;s,t))
= \dim\ker\partial^1 = \left \{ \begin{array}{ll} s+t+1 & \mbox{ if $\car K \mid \gcd(s,t)$}\\ s+t & \mbox{ otherwise}.
\end{array} \right. $$

This completes the proof. \end{proof}

\bigskip

\subsection{$\HH^0(A)$ and $\HH^1(A)$ for the algebras $A=\Gamma(p,q;r)$}

We now turn to the algebras $\Gamma(p,q;r)=K\Q_{(p,q)}/I_r.$ Set $\eta_i=\alpha_i\cdots \alpha_p\delta\alpha_1\cdots
\alpha_{i-1}$ for $1\ppq i\ppq p$ (so that $\eta_1=\gamma\delta$) and $\theta_j=\beta_j\cdots
\beta_q\gamma\beta_1\cdots \beta_{j-1}$ for $1\ppq j\ppq q$ (so that $\theta_1=\delta\gamma$).

\bprop\label{prop:centre of gamma} Let $p\ppq q$ be positive integers and consider the algebra $\Gamma(p,q;r).$ Set
$z=\sum_{i=1}^p\eta_i+\sum_{j=1}^q\theta_j.$ Then \[ \dim \HH^0(\Gamma(p,q;r))= \begin{cases} p+q+r-1&\text{if
$p>1$}\\ q+r+1&\text{if $p=1$ and $q>1$}\\ r+3&\text{if $p=1=q$,} \end{cases}
 \] and a basis for $\HH^0(\Gamma(p,q;r))$ is given by
\begin{align*} &\set{1;(\gamma\delta)^r;z^k,1\ppq k\ppq r-1;\eta_i^r;\theta_j^r\text{ for } 2\ppq i\ppq p \text{ and }
2\ppq j\ppq q}\text{ if $p>1$}\\ &\set{1;(\gamma\delta)^r;z^k,1\ppq k\ppq
r-1;(\delta\gamma)^{r-1}\delta;\theta_j^r\text{ for } 2\ppq j\ppq q}\text{ if $p=1$ and $q>1$}\\
&\set{1;(\gamma\delta)^r;z^k,1\ppq k\ppq r-1;(\gamma\delta)^{r-1}\gamma;(\delta\gamma)^{r-1}\delta}\text{ if $p=1=q$.}
\end{align*}
\eprop

\begin{proof} It is clear that $(\gamma\delta)^r$, $\eta_i^r$ and $\theta_j^r$ are in the centre of
$\Gamma(p,q;r)$ since they are socle elements, and it is easy to check that $(\gamma\delta)^{r-1}\gamma$ and
$(\delta\gamma)^{r-1}\delta$ are in the centre in the appropriate cases.

Conversely, a central element $\zeta$ must be in $\bigoplus_{v=1}^{p+q-1}e_v{\Gamma(p,q;r)}e_v$ and therefore is a linear
combination of  $(\gamma\delta)^k$ for $0\ppq k\ppq r$,   $(\delta\gamma)^k$ for $1\ppq k\ppq r-1$,
$(\gamma\delta)^k\gamma$ and $(\delta\gamma)^k\delta$ for $0\ppq k\ppq r-1$, and $\eta_i^k$ and $\theta_j^k$ for $2\ppq
i\ppq p,$ $2\ppq j\ppq q$ and $1\ppq k\ppq r.$ Writing the equations $\alpha_i\zeta=\zeta\alpha_i$ and
$\beta_j\zeta=\zeta\beta_j$ gives the result (noting that $\eta_i^k=\alpha_i\cdots
\alpha_p(\delta\gamma)^{k-1}\delta\alpha_1\cdots \alpha_{i-1}$ and $\theta_j^k=\beta_j\cdots
\beta_q(\gamma\delta)^{k-1}\gamma\beta_1\cdots \beta_{j-1}$). \end{proof}

We now use the same method as for $\Lambda(p,q;s,t)$ to compute $\HH^1(\Gamma(p,q;r))$, starting with a minimal set $g^2$ of
generators of $I_r$.

\bprop\label{prop:min gens for gammaa} Consider the algebra $\Gamma(p,q;r)$ with $1 \ppq p \ppq q$.

If $p \pgq 2$ then the following elements form a minimal set of generators for the ideal $I_r$: \[\begin{array}{lcl}
g^2_1 & = & (\gamma\delta)^r-(\delta\gamma)^r\\ g^2_i & = & \eta_i^r\alpha_i  \text{ for all $2\ppq i\ppq p-1$}\\ g^2_p
& = & \alpha_p\alpha_1\\ g^2_{p+1} & = & \beta_q\beta_1\\ g^2_{p+j} & = & \theta_j^r \beta_j \text{ for all $2\ppq
j\ppq q-1$}.\\ \end{array}\]

If $p = 1$ then the following elements form a minimal set of generators for the ideal $I_r$: \[\begin{array}{lcl} g^2_0
& = & (\gamma\delta)^r-(\delta\gamma)^r\\ g^2_1 & = & \alpha^2\\ g^2_2 & = & \beta_q\beta_1\\ g^2_{j+1} & = &
\theta_j^r \beta_j  \text{ for all $2\ppq j\ppq q-1$.}\\ \end{array}\] \eprop

The proof of the next result giving the first Hochschild cohomology group is a similar calculation to that of Proposition \ref{prop:HH1 for lambda} and so is left to the reader.

\bprop\label{prop:HH1 for gamma}
If $p \pgq 2$ then \[\dim\HH^1(\Gamma(p,q;r)) = r+1. \] If $p = 1$ and $q\pgq 2$ then
\[\dim\HH^1(\Gamma(1,q;r)) = \left \{ \begin{array}{ll} r+4 & \text{ if $\car K=2$}\\ r+2 & \text{ if $\car K\neq 2$}.
\end{array} \right. \] If $q = 1$ then $p=1$ and \[\dim\HH^1(\Gamma(1,1;r)) = \left \{ \begin{array}{ll} 2r+6 &
\text{ if $\car K =2$}\\ 2r+2 & \text{ if $\car K\neq 2$}. \end{array} \right. \]
\eprop

\bigskip

\subsection{Higher Hochschild cohomology groups for $\Gamma(p,q;r)$}
 In order to distinguish the algebras of the form $\Gamma(p,q;r)$ up to derived equivalence we need the dimensions of the Hochschild cohomology groups up to $\HH^{2p-2}(\Gamma(p,q;r))$. If $p=1$, this is just $\HH^{0}(\Gamma(1,q;r))$ which we already know, so we shall assume that $p>1$ in the remainder of this section. We begin by giving the start of a projective bimodule resolution of $\Gamma(p,q;r)$ to enable us to find these groups. For ease of notation, set $\Gamma = \Gamma(p,q;r)$.

The projective bimodules $P^n$ in a minimal projective bimodule resolution of $\Gamma$ are known from \cite{Happel}; specifically, the multiplicity of $\Gamma e_i\ot e_j\Gamma$ as a direct summand
of $P^n$ is equal to the dimension of $\Ext^n_\Gamma(S_i, S_j)$, where $S_k$ is the simple module at the vertex $k$. We thus define projective $\Gamma$-$\Gamma$-bimodules (equivalently $\Gamma^e$-modules) $P^0, P^1, \ldots, P^{2p}$ which will be the projectives in our minimal projective bimodule resolution for $\Gamma$.

\bd Let $\Gamma = \Gamma(p,q;r)$ with $p>1$. We define projective $\Gamma$-$\Gamma$-bimodules $P^0, P^1, \ldots, P^{2p}$ as follows.
{\allowdisplaybreaks\begin{align*}
&P^0=\bigoplus_{i=1}^{p+q-1}\Gamma e_i\ot e_i\Gamma,\\
&P^1=\bigoplus_{i=1}^{p-1}\Gamma e_i\ot e_{i+1}\Gamma \oplus \Gamma e_p\ot e_1\Gamma \oplus \bigoplus_{j=2}^{q-1}\Gamma e_{p+j-1}\ot e_{p+j}\Gamma \oplus \Gamma e_{p+q-1}\ot e_1\Gamma \oplus \Gamma e_1\ot e_{p+1}\Gamma,\\
&P^{2n}=\bigoplus_{i=2}^{p-n}\Gamma e_i\ot e_{i+n}\Gamma \oplus \bigoplus_{i=p-n+1}^{p-1}\Gamma e_i\ot e_{i+n+1-p}\Gamma \oplus \Gamma e_p\ot e_{n+1}\Gamma \oplus \bigoplus_{j=2}^{q-n}\Gamma e_{p+j-1}\ot e_{p+j+n-1}\Gamma\\&\qquad\qquad\oplus \bigoplus_{j=q-n+1}^{q-1}\Gamma e_{p+j-1}\ot e_{p+j+n-q}\Gamma \oplus \Gamma e_{p+q-1}\ot e_{p+n}\Gamma \oplus \Gamma e_1\ot e_{1}\Gamma,\ \text{ for $1\ppq n<p$,}\\
&P^{2n-1}=\bigoplus_{i=2}^{p-n}\Gamma e_i\ot e_{i+n}\Gamma \oplus \bigoplus_{i=p-n+1}^{p-1}\Gamma e_i\ot e_{i+n-p}\Gamma \oplus \Gamma e_p\ot e_{n}\Gamma\oplus \Gamma e_1\ot e_{n+1}\Gamma \\&\qquad\qquad\oplus \bigoplus_{j=2}^{q-n}\Gamma e_{p+j-1}\ot e_{p+j+n-1}\Gamma\oplus \bigoplus_{j=q-n+2}^{q-1}\Gamma e_{p+j-1}\ot e_{p+j+n-q-1}\Gamma \oplus \Gamma e_{p+q-1}\ot e_{p+n-1}\Gamma \\&\qquad\qquad\oplus \Gamma e_1\ot e_{p+n}\Gamma\oplus \Gamma e_{p+q-n}\ot e_1\Gamma, \text{ for $2\ppq n<p$}\\
&P^{2p-1}=\bigoplus_{i=1}^{p}\Gamma e_i\ot e_{i}\Gamma \oplus \bigoplus_{j=2}^{q-p}\Gamma e_{p+j-1}\ot e_{2p+j-1}\Gamma\oplus \bigoplus_{j=q-p+2}^{q}\Gamma e_{p+j-1}\ot e_{2p+j-q-1}\Gamma \\&\qquad\qquad\oplus \Gamma e_q\ot e_{1}\Gamma\oplus \Gamma e_{1}\ot e_{2p}\Gamma,\\
&P^{2p}=\bigoplus_{i=1}^{p-1}\Gamma e_i\ot e_{i+1}\Gamma \oplus \bigoplus_{j=2}^{q-p}\Gamma e_{p+j-1}\ot e_{2p+j-1}\Gamma\oplus \bigoplus_{j=q-p+1}^{q}\Gamma e_{p+j-1}\ot e_{2p+j-q}\Gamma \\&\qquad\qquad\oplus \Gamma e_p\ot e_{1}\Gamma\oplus \Gamma e_{1}\ot e_{1}\Gamma.
\end{align*}}
\ed

The first maps $d^{i}: P^{i}\rightarrow P^{i-1}$, for $i = 1, 2, 3$ of a minimal projective bimodule resolution are given in \cite{GSn}. We extend the resolution in \cite{GSn} for our algebra $\Gamma$, in a similar way to \cite{ST}, to make the following definition of maps $d^{i}: P^{i}\rightarrow P^{i-1}$ for $i = 1, \ldots , 2p$.

\bd
Let $\Gamma = \Gamma(p,q;r)$ with $p>1$. We define $\Gamma$-$\Gamma$-bimodule homomorphisms $d^{i}: P^{i}\rightarrow P^{i-1}$, for $i = 1, \ldots , 2p$, as follows.

\smallskip

\begin{enumerate}[$\rhd$,itemindent={-1cm}]
\item $\bm{d^{2n-1}:P^{2n-1}\rightarrow P^{2n-2}}$, for $1\ppq n<p$,
{\allowdisplaybreaks\begin{align*}
\bullet \  & e_i\ot e_{i+n}\mapsto e_i\ot \alpha_{i+n-1}-\alpha_i\ot e_{i+n}\quad\text{for $2\ppq i\ppq p-n$}\\
 \bullet \  &e_i\ot e_{i+n-p}\mapsto\\ &\quad(-1)^{i+p+1}\left[ e_i\ot \eta_{i+n-p}^r+(-1)^{(p-i)(n-1)}\alpha_i\cdots \alpha_p\ot \alpha_1\cdots \alpha_{n-p+i-1}-\eta_i^r\ot e_{i+n-p}\right.\\ &\quad -(-1)^n(-1)^{(p-i)(n-1)}\dsty\sum_{m=1}^{p-i}(-1)^{m(n-1)}\alpha_{i}\cdots\alpha_{p-m}\ot \alpha_{n-m+1}\cdots \alpha_p\delta(\gamma\delta)^{r-1}\alpha_1\cdots \alpha_{n-p+i-1}    \\ &\quad\left. +(-1)^n(-1)^{(p-i)(n-1)}\dsty\sum_{m=p-i+2}^{n-1}(-1)^{m(n-1)}\alpha_{i}\cdots\alpha_p(\delta\gamma)^{r-1}\delta\alpha_1\cdots\alpha_{p-m}\ot \alpha_{n-m+1}\cdots  \alpha_{n-p+i-1}    \right] \\&\text{for $p-n+2\ppq i\ppq p$}\\
\bullet \  &e_{p-n+1}\ot e_1\mapsto\\  &\quad e_{p-n+1}\ot \alpha_p+(-1)^n\dsty\sum_{m=1}^{n-1}(-1)^{m(n-1)}\alpha_{p-n+1}\cdots \alpha_{p-m}\ot \alpha_{n-m+1}\cdots \alpha_p\delta(\gamma\delta)^{r-1}-\alpha_{p-n+1}\cdots\alpha_p\ot e_1\\
\bullet \  &e_{1}\ot e_{n+1}\mapsto\\ &\quad (-1)^{n-1}\left[e_{1}\ot \alpha_1\cdots \alpha_n+(-1)^n\dsty\sum_{m=1}^{n-1}(-1)^{m(n-1)}\delta(\gamma\delta)^{r-1}\alpha_{1}\cdots \alpha_{p-m}\ot \alpha_{n-m+1}\cdots \alpha_n+(-1)^n\alpha_{1}\ot e_{n+1}\right]\\
 \bullet \  &e_{p+j-1}\ot e_{p+j+n-1}\mapsto -e_{p+j-1}\ot \beta_{j+n-1}+\beta_j\ot e_{p+j+n-1}\quad\text{for $2\ppq j\ppq q-n$}\\
\bullet \  & e_{p+j-1}\ot e_{p+j+n-q-1}\mapsto\\ &\quad (-1)^{j+q}\left[ e_{p+j-1}\ot \theta_{j+n-q}^r+(-1)^{(q-j)(n-1)}\beta_j\cdots \beta_q\ot \beta_1\cdots \beta_{n+j-q-1}-\theta_j^r\ot e_{p+j+n-q-1}\right.\\ &\quad -(-1)^n(-1)^{(q-j)(n-1)}\dsty\sum_{m=1}^{q-j}(-1)^{m(n-1)}\beta_{j}\cdots\beta_{q-m}\ot \beta_{n-m+1}\cdots \beta_q(\gamma\delta)^{r-1}\gamma\beta_1\cdots \beta_{n-q+j-1}    \\ &\quad\left. +(-1)^n(-1)^{(q-j)(n-1)}\dsty\sum_{m=q-j+2}^{n-1}(-1)^{m(n-1)}\beta_{j}\cdots\beta_q\gamma(\delta\gamma)^{r-1}\beta_1\cdots\beta_{q-m}\ot \beta_{n-m+1}\cdots  \beta_{n-q+j-1}    \right] \\&\text{for $q-n+2\ppq j\ppq q$}\\
\bullet \  &e_{p+q-n}\ot e_1\mapsto\\ &\quad -\left[e_{p+q-n}\ot \beta_q+(-1)^n\dsty\sum_{m=1}^{n-1}(-1)^{m(n-1)}\beta_{q-n+1}\cdots \beta_{q-m}\ot \beta_{n-m+1}\cdots \beta_q(\gamma\delta)^{r-1}\gamma-\beta_{q-n+1}\cdots\beta_q\ot e_1\right]\\
\bullet \  &e_{1}\ot e_{p+n}\mapsto\\ &\quad (-1)^{n}\left[e_{1}\ot \beta_1\cdots \beta_n+(-1)^n\dsty\sum_{m=1}^{n-1}(-1)^{m(n-1)}(\gamma\delta)^{r-1}\gamma\beta_{1}\cdots \beta_{q-m}\ot \beta_{n-m+1}\cdots \beta_n+(-1)^n\beta_{1}\ot e_{p+n}\right]
\end{align*}}

\item $\bm{d^{2n}:P^{2n}\rightarrow P^{2n-1}}$, for $1\ppq n<p$,
{\allowdisplaybreaks \begin{align*}
\bullet \ & e_i\ot e_{i+n}\mapsto\\  &\quad \dsty\sum_{k=0}^r\eta_i^k\ot \eta_{i+n}^{r-k}+\dsty\sum_{k=0}^{r-1}\left[\dsty\sum_{m=0}^{p-i-n}\eta_i^k\alpha_i\cdots \alpha_{i+m}\ot \alpha_{i+m+n+1}\cdots \alpha_p\delta\alpha_1\cdots \alpha_{i+n-1}\eta_{i+n}^{r-k-1}\right.\\  &\quad\left.+\dsty\sum_{m=0}^{i-2}\eta_i^k\alpha_i\cdots \alpha_{p}\delta\alpha_1\cdots\alpha_m\ot \alpha_{m+n+1}\cdots \alpha_{i+n-1}\eta_{i+n}^{r-k-1}\right.\\  &\quad\left.+(-1)^n\dsty\sum_{j=0}^{q-n}\eta_i^k\alpha_i\cdots\alpha_p\beta_1\cdots\beta_j\ot \beta_{j+n+1}\cdots \beta_q\alpha_1\cdots\alpha_{i+n-1}\eta_{i+n}^{r-k-1}\right]\\ & \text{for $2\ppq i\ppq p-n$}\\
\bullet \ & e_i\ot e_{i+n+1-p}\mapsto  e_i\ot \alpha_{i+n-p}-(-1)^n\alpha_i\ot  e_{i+n+1-p} \text{ for $p+1-n\ppq i\ppq p$}\\
\bullet \ & e_1\ot e_1\mapsto\\  &\quad \dsty\sum_{k=0}^{r-1}\left[\dsty\sum_{i=0}^{p-n}\left(\delta(\gamma\delta)^{k}\alpha_1\cdots \alpha_i\ot \alpha_{i+n+1}\cdots \alpha_p(\delta\gamma)^{r-k-1}+(-1)^n(\gamma\delta)^{k}\alpha_1\cdots \alpha_i\ot \alpha_{i+n+1}\cdots \alpha_p(\delta\gamma)^{r-k-1}\delta\right)\right.\\  &\quad\left. +\dsty\sum_{j=0}^{q-n}\left(\gamma(\delta\gamma)^{k}\beta_1\cdots \beta_j\ot \beta_{j+n+1}\cdots \beta_q(\gamma\delta)^{r-k-1}+(-1)^n(\delta\gamma)^{k}\beta_1\cdots \beta_j\ot \beta_{j+n+1}\cdots \beta_q(\gamma\delta)^{r-k-1}\gamma\right) \right]\\
\bullet \ & e_{p+j-1}\ot e_{p+j+n-q}\mapsto  e_{p+j-1}\ot \beta_{j+n-q}-(-1)^n\beta_j\ot e_{p+j+n-q} \text{ for $q+1-n\ppq j\ppq q$}\\
\bullet \ & e_{p+j-1}\ot e_{p+j+n-1}\mapsto \dsty\sum_{k=0}^r\theta_j^k\ot \theta_{j+n}^{r-k}+\dsty\sum_{k=0}^{r-1}\left[\dsty\sum_{m=0}^{q-j-n}\theta_j^k\beta_j\cdots \beta_{j+m}\ot \beta_{j+m+n+1}\cdots \beta_q\gamma\beta_1\cdots \beta_{j+n-1}\theta_{j+n}^{r-k-1}\right.\\ &\quad\left.+\dsty\sum_{m=0}^{j-2}\theta_j^k\beta_j\cdots \beta_{q}\gamma\beta_1\cdots\beta_m\ot \beta_{m+n+1}\cdots \beta_{j+n-1}\theta_{j+n}^{r-k-1}\right.\\ &\quad\left.+(-1)^n\dsty\sum_{i=0}^{p-n}\theta_j^k\beta_j\cdots\beta_q\alpha_1\cdots\alpha_i\ot \alpha_{i+n+1}\cdots \alpha_p\beta_1\cdots\beta_{j+n-1}\theta_{j+n}^{r-k-1}\right]\\ & \text{for $2\ppq j\ppq q-n$}
\end{align*}}

\item $\bm{d^{2p-1}:P^{2p-1}\rightarrow P^{2p-2}}$
{\allowdisplaybreaks \begin{align*}
\bullet \  & e_i\ot e_{i}\mapsto \\  &\quad (-1)^{i}\left[ e_i\ot \eta_{i}^r+(-1)^{(p-i)(p-1)}\alpha_i\cdots \alpha_p\ot \alpha_1\cdots \alpha_{i-1}-\eta_i^r\ot e_{i}\right.\\  &\quad -(-1)^p(-1)^{(p-i)(p-1)}\dsty\sum_{m=1}^{p-i}(-1)^{m(p-1)}\alpha_{i}\cdots\alpha_{p-m}\ot \alpha_{p-m+1}\cdots \alpha_p\delta(\gamma\delta)^{r-1}\alpha_1\cdots \alpha_{i-1}    \\  &\quad\left. +(-1)^p(-1)^{(p-i)(p-1)}\dsty\sum_{m=p-i+2}^{p-1}(-1)^{m(p-1)}\alpha_{i}\cdots\alpha_p(\delta\gamma)^{r-1}\delta\alpha_1\cdots\alpha_{p-m}\ot \alpha_{p-m+1}\cdots  \alpha_{i-1}    \right] \\&\text{for $2\ppq i\ppq p$}\\
\bullet \  & e_{1}\ot e_1\mapsto \\  &\quad e_{1}\ot \gamma-\dsty\sum_{m=1}^{p-1}(-1)^{m(p-1)}\alpha_{1}\cdots \alpha_{p-m}\ot \alpha_{p-m+1}\cdots \alpha_p\delta(\gamma\delta)^{r-1}\\ & \qquad+(-1)^p\sum_{m=1}^{p-1}(-1)^{m(p-1)}\delta(\gamma\delta)^{r-1}\alpha_{1}\cdots \alpha_{p-m}\ot \alpha_{p-m+1}\cdots \alpha_p+(-1)^p\gamma\ot e_1\\
\bullet \  &  e_{p+j-1}\ot e_{2p+j-1}\mapsto -e_{p+j-1}\ot \beta_{p+j-1}+\beta_j\ot e_{2p+j-1}\quad\text{for $2\ppq j\ppq q-p$}\displaybreak\\
\bullet \  &  e_{p+j-1}\ot e_{2p+j-q-1}\mapsto \\  &\quad (-1)^{j+q}\left[ e_{p+j-1}\ot \theta_{p+j-q}^r+(-1)^{(q-j)(p-1)}\beta_j\cdots \beta_q\ot \beta_1\cdots \beta_{p+j-q-1}-\theta_j^r\ot e_{2p+j-q-1}\right.\\ & \qquad -(-1)^p(-1)^{(q-j)(p-1)}\dsty\sum_{m=1}^{q-j}(-1)^{m(p-1)}\beta_{j}\cdots\beta_{q-m}\ot \beta_{p-m+1}\cdots \beta_q(\gamma\delta)^{r-1}\gamma\beta_1\cdots \beta_{p-q+j-1}    \\  &\quad \left. +(-1)^p(-1)^{(q-j)(p-1)}\dsty\sum_{m=q-j+2}^{p-1}(-1)^{m(p-1)}\beta_{j}\cdots\beta_q\gamma(\delta\gamma)^{r-1}\beta_1\cdots\beta_{q-m}\ot \beta_{p-m+1}\cdots  \beta_{p-q+j-1}    \right] \\&\text{for $q-p+2\ppq j\ppq q$}\\
\bullet \  & e_{q}\ot e_1\mapsto \\ &\quad  -\left[e_{q}\ot \beta_q+(-1)^p\dsty\sum_{m=1}^{p-1}(-1)^{m(p-1)}\beta_{q-p+1}\cdots \beta_{q-m}\ot \beta_{p-m+1}\cdots \beta_q(\gamma\delta)^{r-1}\gamma-\beta_{q-p+1}\cdots\beta_q\ot e_1\right]\\
\bullet \  & e_{1}\ot e_{2p}\mapsto \\  &\quad  (-1)^{p}\left[e_{1}\ot \beta_1\cdots \beta_p+(-1)^p\dsty\sum_{m=1}^{p-1}(-1)^{m(p-1)}(\gamma\delta)^{r-1}\gamma\beta_{1}\cdots \beta_{q-m}\ot \beta_{p-m+1}\cdots \beta_p+(-1)^p\beta_{1}\ot e_{2p}\right]\end{align*}}

\item $\bm{d^{2p}:P^{2p}\rightarrow P^{2p-1}}$
{\allowdisplaybreaks\begin{align*}
\bullet\ & e_i\ot e_{i+1}\mapsto  e_i\ot \alpha_{i}-(-1)^p\alpha_i\ot e_{i+1}\qquad\text{for $1\ppq i\ppq p-1$}\\
\bullet\ & e_p\ot e_1\mapsto  e_p\ot \alpha_p-(-1)^p\alpha_p\ot e_1\\
\bullet\ & e_1\ot e_1\mapsto\\  &\quad\dsty\sum_{k=0}^{r-1}\left[(-1)^p\delta(\gamma\delta)^k\ot(\delta\gamma)^{r-k-1}+(\gamma\delta)^k\ot (\delta\gamma)^{r-k-1}\delta\right.\\  &\quad -\left.\dsty\sum_{j=0}^{q-p}\left(\gamma(\delta\gamma)^{k}\beta_1\cdots \beta_j\ot \beta_{j+p+1}\cdots \beta_q(\gamma\delta)^{r-k-1}+(-1)^p(\delta\gamma)^{k}\beta_1\cdots \beta_j\ot \beta_{j+p+1}\cdots \beta_q(\gamma\delta)^{r-k-1}\gamma\right) \right]\\
\bullet\ & e_{p+j-1}\ot e_{p+j+p-q}\mapsto  e_{p+j-1}\ot \beta_{j+p-q}-(-1)^p\beta_j\ot e_{p+j+p-q}\qquad\text{for $q+1-p\ppq j\ppq q$}\\
\bullet\ & e_{p+j-1}\ot e_{2p+j-1}\mapsto\\  &\quad \dsty\sum_{k=0}^r\theta_j^k\ot \theta_{j+p}^{r-k}+\dsty\sum_{k=0}^{r-1}\left[\dsty\sum_{m=0}^{q-j-p}\theta_j^k\beta_j\cdots \beta_{j+m}\ot \beta_{j+m+p+1}\cdots \beta_q\gamma\beta_1\cdots \beta_{j+p-1}\theta_{j+p}^{r-k-1}\right.\\&\quad\left.+\dsty\sum_{m=0}^{j-2}\theta_j^k\beta_j\cdots \beta_{q}\gamma\beta_1\cdots\beta_m\ot \beta_{m+p+1}\cdots \beta_{j+p-1}\theta_{j+p}^{r-k-1}-\theta_j^k\beta_j\cdots\beta_q\ot\beta_1\cdots\beta_{j+p-1}\theta_{j+p}^{r-k-1}\right]\\ &\quad \text{for $2\ppq j\ppq q-p$}.
\end{align*}}
\end{enumerate}
\ed

It remains to show that the projective bimodules and homomorphisms that we have defined do indeed give the start of a minimal projective bimodule resolution of $\Gamma$.

\bt
With the above notation,
$$\xymatrix{
\cdots\ar[r] & P^{2p}\ar[r]^{d^{2p}} & P^{2p-1}\ar[r]^>>>>{d^{2p-1}} & & \cdots\ar[r] & P^1\ar[r]^{d^1} & P^0\ar[r] & \Gamma\ar[r] & 0
}$$
is the beginning of a minimal projective resolution of $\Gamma$ as a $\Gamma$-$\Gamma$-bimodule (when $p>1$).
\et

\begin{proof}
It may be verified directly from the definitions that $d^2 = 0$ and thus we have a complex.
The strategy for proving exactness is identical to that of \cite[Theorem 1.6]{ST} (and see \cite[Proposition 2.8]{GSn}), whereby we apply $(\Gamma/\rad(\Gamma)\otimes -)$ to the complex and show that the resulting sequence is a minimal projective resolution of $\Gamma/\rad(\Gamma)$ as a right $\Gamma$-module. Minimality is then immediate since we know that the projectives are those of a minimal projective resolution of $\Gamma$ as a $\Gamma$-$\Gamma$-bimodule from \cite{Happel}.
\end{proof}

We are now in a position to give the dimensions of the Hochschild cohomology groups up to $\HH^{2p-2}(\Gamma)$. We only give those in even degree since we shall not need the others. The details of the proof are left to the reader.

\bt\label{thm:higher HH for gamma}
For $2\ppq n<p\ppq q$ we have
\[ \dim \HH^{2n-2}(\Gamma)=
\begin{cases}
r&\text{ if $n$ is odd and $\car K\nmid 2r$}\\
r+1&\text{ if $n$ is odd and $\car K\mid 2r$}\\
r&\text{ if $n$ is even and $\car K\neq 2$}\\
r+1&\text{ if $n$ is even and $\car K= 2$,}\\
\end{cases}
 \] and for $2\ppq p<q$ we have
\[ \dim \HH^{2p-2}(\Gamma)=
\begin{cases}
r-1&\text{ if $p$ is odd and $\car K\nmid 2r$}\\
r&\text{ if $p$ is odd, $\car K\neq 2$ and $\car K\mid 2r$}\\
r+1&\text{ if $p$ is even and $\car K\neq  2$}\\
r+2&\text{ if  $\car K= 2$}.
\end{cases}
 \]
\et

\bigskip

\section{Derived equivalence classes}\label{sec:final}

 It was shown in \cite{BHS} that two algebras of the form $A(p,q) = \Gamma(p,q;1)$ or $\Lambda(n) = \Lambda(1,n;2,2)$ are derived equivalent if and only if they are isomorphic. The main result in this section is to extend this  to all algebras of the form $\Gamma(p,q;r)$ and $\Lambda(p,q;s,t)$, and hence to all basic indecomposable finite-dimensional $K$-algebras $A$ which are symmetric special biserial algebras with at most one non-uniserial indecomposable projective module.

We start with some properties of these algebras, all of which are invariants under derived equivalence.

\bprop\label{prop:properties} Suppose $1 \ppq p \ppq q$. The algebras $\Gamma(p,q;r)$ and $\Lambda(p,q;s,t)$ have the following properties.
\begin{enumerate}[\bf (1)]
\item The number of simples of
  $\begin{cases}
  \Gamma(p,q;r) \mbox{ is } p+q-1,\\
  \Lambda(p,q;s,t) \mbox{ is } p+q-1.
  \end{cases}$
\item The Cartan invariants of
  $\begin{cases}
  \Gamma(p,q;r) \mbox{ are }
     \begin{cases}
     \underbrace{1, 1, \dots , 1,}_{p+q-2} 4r & \text{if $r(p+q-2)$ is even,}\\
     \underbrace{1, 1, \dots , 1,}_{p+q-3} 2, 2r & \text{if $r(p+q-2)$ is odd,}\\
     \end{cases}
   \\
  \Lambda(p,q;s,t) \mbox{ are } \underbrace{1, 1, \dots , 1,}_{p+q-2} s+t+(p+q-2)st.
  \end{cases}$
\item The Cartan determinant of
  $\begin{cases}
  \Gamma(p,q;r) \mbox{ is } 4r,\\
  \Lambda(p,q;s,t) \mbox{ is } s+t+(p+q-2)st.
  \end{cases}$
\end{enumerate}
\eprop

\begin{proof}
\begin{enumerate}[\bf (1)]
\item This is immediate from the number of vertices of the quiver ${\mathcal Q}_{(p,q)}$.
\item Let ${\mathcal I}_m$ be the $m\times m$ identity matrix, let ${\mathcal J}_n$ be the $n \times n$ matrix with all entries equal to 1, and set $u = p+q-2$. We start with the algebra $\Gamma(p,q;r)$. The Cartan matrix of $\Gamma(p,q;r)$ is the $(p+q-1)\times (p+q-1)$ matrix
$$C_\Gamma = \left [ \begin{array}{cc}
4r & 2r \ \ \ \ \ \ \dots \ \ \ \ \ \ 2r\\
2r & \\
\vdots & {\mathcal I}_{u}+r{\mathcal J}_{u}\\
2r & \\
\end{array}\right ].$$
The Smith normal form for $C_\Gamma$ is
$$\begin{cases}
  \left [ \begin{array}{cc}
  {\mathcal I}_{u} & 0\\
  0 & 4r\\
  \end{array}\right ] & \text{if $ru$ is even} \\
  \ \\
  \left [ \begin{array}{ccc}
  {\mathcal I}_{u-1} & 0 & 0\\
  0 & 2 & 0\\
  0 & 0 & 2r\\
  \end{array}\right ] & \text{if $ru$ is odd,}
  \end{cases}$$
and thus the Cartan invariants of $\Gamma(p,q;r)$ are
$\begin{cases}
\underbrace{1, 1, \dots , 1,}_{p+q-2} 4r & \text{if $ru$ is even}\\
\underbrace{1, 1, \dots , 1,}_{p+q-3} 2, 2r & \text{if $ru$ is odd}.\\
\end{cases}$
\\

Now consider the algebra $\Lambda(p,q;s,t)$. The Cartan matrix of $\Lambda(p,q;s,t)$ is the $(p+q-1)\times (p+q-1)$ matrix
$$C_\Lambda = \left [ \begin{array}{ccc}
s+t & t \ \ \ \ \ \ \dots \ \ \ \ \ \ t & s \ \ \ \ \ \ \dots \ \ \ \ \ \ s\\
t & & \\
\vdots & {\mathcal I}_{p-1}+t{\mathcal J}_{p-1} & 0\\
t & &\\
s & &\\
\vdots & 0 & {\mathcal I}_{q-1}+s{\mathcal J}_{q-1}\\
s & &\\
\end{array}\right ].$$
The Smith normal form for $C_\Lambda$ is
$$\left [ \begin{array}{cc}
{\mathcal I}_{u} & 0\\
0 & s+t+ust\\
\end{array}\right ],$$
so the Cartan invariants of $\Lambda(p,q;s,t)$ are $\underbrace{1, 1, \dots , 1,}_{p+q-2} s+t+(p+q-2)st.$
\item This is immediate from (2).
\qedhere\end{enumerate}
\end{proof}

We now consider isomorphism classes of algebras of the same form. It is clear that $\Lambda(q,q;s,t) = \Lambda(q,q;t,s)$, and the next result shows that, with this one exception,  two algebras both of the form $\Gamma(p,q;r)$ or of the form $\Lambda(p,q;s,t)$ are pairwise non-isomorphic.

\bt \begin{enumerate}[\bf (1)]
\item The algebras of the form $\Gamma(p,q;r)$ (with $1 \ppq p \ppq q$) are pairwise non-isomorphic.
\item The algebras of the form $\Lambda(p,q;s,t)$ (with $1 \ppq p \ppq q$) are pairwise non-isomorphic with the exception that $\Lambda(q,q;s,t) = \Lambda(q,q;t,s)$.
\end{enumerate}
\et

\begin{proof} \begin{enumerate}[\bf (1)]
\item First, suppose that the algebras $\Gamma(p,q;r)$ and $\Gamma(p',q';r')$ are isomorphic with $1 \ppq p \ppq q$, $1 \ppq p' \ppq q'$. Since both algebras are basic, the quivers are uniquely determined and hence ${\mathcal Q}_{(p,q)} = {\mathcal Q}_{(p',q')}$. Thus $p = p'$ and $q = q'$. From the Cartan determinant in Proposition~\ref{prop:properties}{\bf (3)}, we have that $r = r'$.
\item Suppose that $\Lambda(p,q;s,t) \cong \Lambda(p',q';s',t')$ with $1 \ppq p \ppq q$, $1 \ppq p' \ppq q'$. Since both algebras are basic, we again have that $p = p'$ and $q = q'$. Then, using the zero-th Hochschild cohomology group from Proposition~\ref{prop:centre of lambda}, we have $s+t = s'+t'$. Equality of the Cartan invariants from Proposition~\ref{prop:properties}{\bf (2)} gives that $st = s't'$. Hence $s', t'$ are the two roots of the equation $x^2-(s+t)x+st = 0$. Thus we have either $s=s'$ and $t=t'$ (which gives us the algebra $\Lambda(p,q;s,s)$), or $s=t'$ and $t=s'$. In the latter case we have the algebras $\Lambda(p,q;s,t)$ and $\Lambda(p,q;t,s)$, and it remains to show $p=q$ when $s \neq t$. The $K$-dimension of $\Lambda(p,q;s,t)$ is $tp^2+sq^2+p+q-2$ and that of $\Lambda(p,q;t,s)$ is $sp^2+tq^2+p+q-2$. Thus $tp^2+sq^2 = sp^2+tq^2$ so that $(t-s)(p+q)(p-q) = 0$. Since $p+q >0$ and $s \neq t$ we have $p=q$, which is precisely the case $\Lambda(q,q;s,t) = \Lambda(q,q;t,s)$.
\qedhere\end{enumerate}
\end{proof}

Our next theorem classifies up to derived equivalence, all basic indecomposable finite-dimensional $K$-algebras $A$ which are symmetric special biserial algebras with at most one non-uniserial indecomposable projective module.

\bt\label{thm:derived equiv}\begin{enumerate}[\bf (1)]
\item An algebra of the form $\Lambda(p,q;s,t)$ (with $1 \ppq p \ppq q$) is derived equivalent to exactly one
algebra in the following list:
\begin{enumerate}[\bf(a)]
\item $\Lambda(1,p+q-1;s,t)$ with $2\ppq s\ppq t$,
\item $N_{M}^{p+q-1}$ with $p+q > 2$ and $\min(s,t)=1$, $\max(s,t)=M$.
\end{enumerate}
\item An algebra of the form $\Gamma(p,q;r)$ (with $1 \ppq p \ppq q$) is derived equivalent to an algebra of the
 form $\Lambda(p,q;s,t)$ if and only if they are isomorphic.
 This is only the case for $\Gamma(1,1;1)\cong \Lambda(1,1;2,2)$ and $\car K \neq 2$.
\item The algebras $\Gamma(p,q;r)$ and $\Gamma(p',q';r')$ (with $1 \ppq p\ppq q$ and $1 \ppq p'\ppq q'$) are derived equivalent if and only if $(p,q,r)=(p',q',r').$
\end{enumerate} \et

\begin{proof}
\begin{enumerate}[\bf (1)]
\item The algebra $\Lambda(p,q;s,t)$ is the generalised Brauer tree algebra associated to the Brauer tree in Figure 1:
\begin{center}
\begin{tabular}{ccc}
{$ \xygraph{
!{<0cm,0cm>;<1.7cm,0cm>:<0cm,1.7cm>::}
!{(0,0) }*+{{}_{a}}="a"
!{(1,0) }*+{{}_{b}}="b"
!{(0,-1) }*+{ }="vx"
!{(-0.5,-0.87) }*+{ }="vxx"
!{(-0.5,0.87) }*+{ }="vxp"
!{(1.5,0.87) }*+{ }="vxpp"
!{(1.5,-0.87) }*+{ }="vxpqq"
!{(1,-1) }*+{}="vxpq"
!{(-0.25,0.43) }*+{ }="A"
!{(-0.25,-0.43) }*+{ }="B"
!{(1.25,0.43) }*+{ }="C"
!{(1.25,-0.43) }*+{ }="D"
"a"-"b" ^{1}
"a"-"vx" _{2}
"a"-"vxx" _{3}
"a"-"vxp" ^{p}
"b"-"vxpp" ^{p+1}
"b"-"vxpqq" ^{p+q-2}
"b"-"vxpq" _{p+q-1}
"A":@{.}@/_/"B"
"D":@{.}@/_/"C"
} $}
&\hspace*{2cm}&
{$ \xygraph{
!{<0cm,0cm>;<1.7cm,0cm>:<0cm,1.7cm>::}
!{(0,0) }*+{{}_{c}}="c"
!{(1,0) }*+{{}_{d}}="d"
!{(1,1) }*+{ }="vx"
!{(1.5,0.87) }*+{ }="vxpp"
!{(1.5,-0.87) }*+{ }="vxpqq"
!{(1,-1) }*+{}="vxpq"
!{(-0.25,0.43) }*+{ }="A"
!{(-0.25,-0.43) }*+{ }="B"
!{(1.25,0.43) }*+{ }="C"
!{(1.25,-0.43) }*+{ }="D"
"c"-"d" ^{1}
"d"-"vx" ^{2}
"d"-"vxpp" ^{3}
"d"-"vxpqq" ^{p+q-2}
"d"-"vxpq" _{p+q-1}
"D":@{.}@/_/"C"
} $}\\\\
Figure 1 && Figure 2
\end{tabular}
\end{center}
in which the vertices $a$ and $b$ have multiplicities respectively $s$ and $t$ (we refer to \cite[Section 4.18]{B} or \cite{MH} for the definition of a Brauer tree algebra and a generalised Brauer tree algebra).

It was proved in \cite[Theorem 9.7]{MH} that generalised Brauer tree algebras up to derived equivalence depend only on the number of edges in the graph and the set of multiplicities. Therefore $\Lambda(p,q;s,t)$ is derived equivalent to the generalised Brauer tree algebra associated to the Brauer tree in Figure 2,
in which the vertices  $c$ and $d$ have multiplicities respectively $m=\min(s,t)$ and $M=\max(s,t)$ and $\set{a,b}=\set{c,d}.$ This algebra is equal to
\begin{itemize}
\item either $K\Delta_{p+q-1}/L_M=N_{M}^{p+q-1}$ if $m=1$, that is, if  $(s,t)=(M,1)$ or $(s,t)=(1,M)$, with $M\pgq 1,$
\item or $\Lambda(1,p+q-1;m,M)$ if $m>1$, that is, if $s\pgq 2$ and $t\pgq 2.$
\end{itemize} Moreover, none of  these algebras are  derived equivalent, again by \cite[Theorem 9.7]{MH}.
\item First, we show that the algebras $\Gamma(1,1;1)$ and $\Lambda(1,1;2,2)$ are isomorphic when $\car K \neq 2$. Since $K$ is algebraically closed, let $\varepsilon$ be a square root of $-1$ in $K$. Then the map
    $$\varphi: \Gamma(1,1;1) \to \Lambda(1,1;2,2) \mbox{ given by }
    \left \{ \begin{array}{l}
    \alpha \mapsto \alpha + \varepsilon \beta\\
    \beta \mapsto \alpha - \varepsilon \beta
    \end{array} \right.$$
is an isomorphism of algebras.

\medskip

Suppose that there is a derived equivalence between the algebras $\Gamma(p,q;r)$ and $\Lambda(p',q';s,t)$. Then the algebras have the same number of simple modules so, from Proposition~\ref{prop:properties}, we have $p+q = p'+q'$.
From {\bf (1)}, the algebra $\Lambda(p',q';s,t)$ is derived equivalent to exactly one algebra in the list above. Moreover, $s+t=m+M$ where $m=\min(s,t)$ and $M=\max(s,t).$

\medskip

\noindent {\bf Case I}: $m=\min(s,t)=1$.

Let $\Gamma = \Gamma(p,q;r)$ and $\Lambda = \Lambda(p',q';s,t)$. From Propositions~\ref{prop:centre of lambda} and \ref{prop:HH1 for lambda}, $\dim\HH^0(\Lambda) = p+q+M-1$ and $\dim\HH^1(\Lambda) = M$.
We first assume that $p=1$, so that $q>1$. From Proposition~\ref{prop:centre of gamma}, $\dim\HH^0(\Gamma) = q+r+1$ so that $M = r+1$. However, Proposition~\ref{prop:HH1 for gamma} gives
$$\dim\HH^1(\Gamma) = \begin{cases}
r+4  & \text{if $\car K = 2$}\\
r+2  & \text{if $\car K \neq 2$}
\end{cases}$$
which is a contradiction. So $\Gamma(1,q;r)$ is not derived equivalent to $\Lambda(p',q';s,t)$.

On the other hand, if $p>1$ then $\dim\HH^0(\Gamma) = p+q+r-1$ so that $M=r$. But $\dim\HH^1(\Gamma) = r+1$ so that $M = r+1$, a contradiction. Again $\Gamma(p,q;r)$ is not derived equivalent to $\Lambda(p',q';s,t)$.

\bigskip

\noindent {\bf Case II}: $m=\min(s,t)>1.$

We use Propositions~\ref{prop:centre of lambda}, \ref{prop:HH1 for lambda}, \ref{prop:centre of gamma} and \ref{prop:HH1 for gamma} without further comment. We begin with the case $p=1$. If $q=1$ then $\dim\HH^0(\Lambda) = s+t$ and $\dim\HH^0(\Gamma) = r+3$. From Proposition~\ref{prop:properties}{\bf (3)}, the Cartan determinant of $\Gamma$ is $4r$ and of $\Lambda$ is $s+t$. Hence $r+3 = s+t = 4r$ so that $r=1$ and $s=2=t$. If $\car K \neq 2$ then we have an isomorphism $\Gamma(1,1;1) \cong \Lambda(1,1;2,2)$ from above. If $\car K = 2$, then $\dim\HH^1(\Lambda) = 5$ and $\dim\HH^1(\Gamma) = 8$, which is a contradiction and $\Gamma(1,1;1)$ is not derived equivalent to $\Lambda(1,1;2,2)$. In the case where $q>1$ (keeping $p = 1$), a similar consideration of the zero-th and first Hochschild cohomology groups and the Cartan determinant shows that $\Gamma(1,q;r)$ is not derived equivalent to $\Lambda(p',q';s,t)$.

Now suppose that $p > 1$. We have $\dim\HH^0(\Lambda) = p+q+s+t-2$ and $\dim\HH^0(\Gamma) = p+q+r-1$ so that $r = s+t-1$. So $\dim\HH^1(\Gamma) = r+1 = s+t$, and we must have $\car K = \ell | \gcd(s,t)$. From Proposition~\ref{prop:properties}{\bf (3)}, the Cartan determinant of $\Gamma$ is $4r$ and  the Cartan determinant of $\Lambda$ is $s+t+(p+q-2)st$ so that $3(s+t) = (p+q-2)st - 4$. Thus $\ell | 4$ so that $\ell =2$, and $s,t$ are both even.
Thus we are in the situation where $\Gamma(p,q;r)$ is derived equivalent to $\Lambda(p',q';s,t)$, $\car K = 2$, both $s, t$ are even, $r = s+t-1$, and $p+q = p'+q'$. 
We shall use K\"ulshammer invariants and the same arguments as in \cite[Subsection 4.5.2]{HZ} for this case. Recall from Proposition \ref{prop:algebrassatisfyproperty} and its proof that the algebras  $\Lambda(p',q';s,t)$ and $\Gamma(p,q;r)$ are symmetric; moreover, once we have fixed a $K$-basis of paths for the socle of a symmetric algebra $A$ and completed it to a $K$-basis $\B_A$ of paths for $A$, the linear map $f:A\rightarrow K$ that is defined on $\B_A$ by sending socle elements to  $1$ and the rest to $0$ defines a symmetric non-degenerate associative bilinear form  on $A$. Orthogonals will be taken with respect to this bilinear form. Let $\kappa(A)$ be the commutator subspace of $A$ and, for any non-negative integer $n$, define $T_n(A)=\set{x\in A;x^{2^n}\in \kappa(A)}.$ It was proved in \cite{Z} that the generalised Reynolds ideals (or K\"ulshammer invariants) $T_n(A)^\perp$ are derived invariant.
Note that $\soc(A)\subseteq T_n(A)^\perp\subsetneq Z(A)$ for every $n$. It is well-known that the centre $Z(A)$ is a derived invariant. Given a vector space $V$, let $\B_V$ denote a basis of $V.$
\begin{itemize}
\item We start with the algebra $\Gamma=\Gamma(p,q;r)$. A basis $\B_\Gamma$ of paths of $\Gamma$ is given by the union over all $i,j$ with $1\ppq i,j\ppq p+q-1$ of all paths from $e_i$ to $e_j$ of length at most $ (p+q)r$ except $(\delta\gamma)^r$.   Recall from Proposition \ref{prop:algebrassatisfyproperty} the basis $\B_{\soc(\Gamma)}=\set{(\gamma\delta)^r;\eta_i^r,2\ppq i\ppq p;\theta_i^r,2\ppq j\ppq q}\subset\B_\Gamma$  of $\soc(\Gamma)$ and from Proposition \ref{prop:centre of gamma} the basis $\B_{Z(\Gamma)}=\set{1;(\gamma\delta)^r;z^k,1\ppq k\ppq r-1;\eta_i^r;\theta_j^r\text{ for } 2\ppq i\ppq p \text{ and }2\ppq j\ppq q}$ of $Z(\Gamma)$ where $z=\sum_{i=1}^p\eta_i+\sum_{j=1}^q\theta_j.$
Then $\dim \kappa(\Gamma)=r\left((p+q)^2-1\right)-1$ and a basis of $\kappa(\Gamma)$ is given by
\begin{align*}
\B_{\kappa(\Gamma)}=\{\gamma(\delta\gamma)^k;\delta(\gamma\delta)^k;&(\gamma\delta)^{k+1}-(\delta\gamma)^{k+1}; \eta_i^{k+1}-\eta_1^{k+1}, 1\ppq i\ppq p;\theta_j^{k+1}-\theta_1^{k+1},2\ppq j\ppq q;\\&\text{ for }0\ppq k\ppq r-1\}\cup\set{b\in \B_\Gamma;\mo(b)\neq \mt(b)}.
\end{align*} Now, as in \cite[Subsection 4.5.2]{HZ}, work in $\Gamma/\kappa(\Gamma)$ to find a basis $\B_{T_1(\Gamma)}=\B_{\kappa(\Gamma)}\cup\set{(\gamma\delta)^k;\frac{r+1}{2}\ppq k\ppq r}$ of $T_1(\Gamma)$ (recall that $r$ is odd), then work in $Z(\Gamma)/\soc(\Gamma)$ to find a basis $\B_{T_1(\Gamma)^\perp}=\B_{\soc(\Gamma)}\cup \set{z^k;\frac{r+1}{2}\ppq k\ppq r-1}$ for $T_1(\Gamma)^\perp$ so that $\Gamma':=Z(\Gamma)/T_1(\Gamma)^\perp$ has basis $\B_{\Gamma'}=\set{1;z^k;1\ppq k\ppq \frac{r-1}{2}}.$ We now consider the Jacobson radical $\rrad_{\Gamma'}$ of the algebra $\Gamma'$, which is spanned by $\set{z^k;1\ppq k\ppq \frac{r-1}{2}},$ and its square $\rrad_{\Gamma'}^2$, which is spanned by $\set{z^k;2\ppq k\ppq \frac{r-1}{2}}$ so that $\dim_K\rrad_{\Gamma'}/\rrad_{\Gamma'}^2=1.$
\item We now turn to the algebra $\Lambda=\Lambda(p',q';s,t)$. Since they are derived equivalent, we may assume that $\Lambda=\Lambda(1,n;m,M)$ with $n=p+q-1$ (to simplify notation). Note that $m$ and $M$ are even. Set $\alpha=\alpha_1.$ We follow the same method as for $\Gamma,$ using Propositions~\ref{prop:algebrassatisfyproperty} and \ref{prop:centre of lambda}.

\begin{align*}
&\B_\Lambda=\bigcup_{1\ppq i,j\ppq p+q-1}\set{\text{all paths  from $e_i$ to $e_j$ of length at most } nM \text{ that do not contain }\alpha }\\&\qquad\qquad\qquad\qquad\qquad\cup\set{\alpha^k,1\ppq k\ppq m-1}\\
&\B_{\soc(\Lambda)}=\set{\delta_i^m, 1\ppq i\ppq n }\\
&\B_{Z(\Lambda)}=\set{1;\alpha^k,1\ppq k\ppq m; y^\ell,1\ppq \ell\ppq t-1;\delta_j^M,2\ppq j\ppq q } \text{ where $y=\sum_{j=1}^q\delta_j$}\\
&\B_{\kappa(\Lambda)}=\set{\delta^\ell-\delta_j^\ell;2\ppq j\ppq q,1\ppq \ell\ppq M}\cup\set{b\in\B_\Lambda;\mo(b)\neq \mt(b)}\\
&\B_{T_1(\Lambda)}=\B_{\kappa(\Lambda)}\cup \set{\alpha^k,\frac{m}{2}<k\ppq m;\delta^\ell,\frac{M}{2}<\ell\ppq M;\alpha^{m/2}+\delta^{M/2}}\\
&\B_{T_1(\Lambda)^\perp}=\B_{\soc(\Lambda)}\cup\set{\alpha^k,\frac{m}{2}<k\ppq m-1;y^\ell,\frac{M}{2}<\ell\ppq M-1;\alpha^{m/2}+y^{M/2}}.
\end{align*}
So the Jacobson radical $\rrad_{\Lambda'}$ of  $\Lambda':=Z(\Lambda)/T_1(\Lambda)^\perp$ has basis $\set{\alpha^k,1\ppq k\ppq \frac{m}{2};y^\ell,1\ppq \ell<\frac{M}{2}}$ and  $\rrad_{\Lambda'}^2$ has basis $\set{\alpha^k,2\ppq k\ppq \frac{m}{2};y^\ell,2\ppq \ell<\frac{M}{2}}$ so that \[\dim_K \rrad_{\Lambda'}/\rrad_{\Lambda'}^2=
\begin{cases}
2\text{ if $M\pgq 4$}\\1\text{ if $M=2$ (and therefore $m=2$)}.
\end{cases}
\]
\end{itemize} Since we assumed that $\Lambda$ and $\Gamma$ are derived equivalent, the algebras $\Lambda'$ and $\Gamma'$ are isomorphic and hence we have $\dim_K\rrad_{\Gamma'}/\rrad_{\Gamma'}^2=\dim_K \rrad_{\Lambda'}/\rrad_{\Lambda'}^2$, which implies that $m=2=M$, that is, $s=2=t$.  Therefore  $r=s+t-1=3$. We now use the Cartan determinants: $0=\det C_\Gamma-\det C_\Lambda=4r-(s+t+(p+q-2)st)=12-(4+(p+q-2)4)=4(4-(p+q))$ so that  $p+q=4.$ Since  $1<p\ppq q$, we must have $p=2=q.$ Therefore $\Gamma=\Gamma(2,2;3)$ and $\Lambda$ is derived equivalent to $\Lambda(2,2;2,2).$ However, it was shown in \cite[Section 3]{H} that $\Gamma(2,2;3)=\D(3\A)_1^3$ and $\Lambda(2,2;2,2)=\D(3\A)_2^{2,2}$ are not derived equivalent.
 Therefore $\Lambda$ and $\Gamma$ are not derived equivalent.

\item If the algebras $\Gamma(p,q;r)$ and $\Gamma(p',q';r')$ (with $p\ppq q$ and $p'\ppq q'$) are derived equivalent, then from the Cartan determinant and number of simples (Proposition~\ref{prop:properties}) we know that $r=r'$ and $p+q=p'+q'$. Assume for contradiction that $(p,q)\neq (p',q')$. We may suppose that $p<p'.$ It follows that $p<q$ (for otherwise $p=q$ and hence
    $q'=(p'+q')-p'=(p+q)-p'=2p-p'<2p'-p'=p'$, a contradiction). Using Theorem~\ref{thm:higher HH for gamma} and Proposition~\ref{prop:centre of gamma}, we then have \[ \HH^{2p-2}\left(\Gamma(p,q;r)\right)\neq \HH^{2p-2}\left(\Gamma(p',q';r')\right) \] which contradicts the fact that the algebras are derived equivalent. Thus $(p,q) = (p',q')$.
\qedhere\end{enumerate}
\end{proof}

\section{Stable equivalence of Morita type classes}\label{sec:extrafinal}

Finally, we  give a classification up to stable equivalence of Morita type of  all algebras of the
form $\Gamma(p,q;r)$ and $\Lambda(p,q;s,t)$, and hence of all basic indecomposable
finite-dimensional $K$-algebras $A$ which are symmetric special biserial algebras with at most one
non-uniserial indecomposable projective module. This is based on the classification up to derived
equivalence of Theorem \ref{thm:derived equiv}.

\bt\label{thmSEMT}\begin{enumerate}[\bf (1)]
\item An algebra of the form $\Lambda(p,q;s,t)$ (with $1 \ppq p \ppq q$) is stably equivalent of
  Morita type to exactly one
algebra in the following list:
\begin{enumerate}[\bf(a)]
\item $\Lambda(1,p+q-1;s,t)$ with $2\ppq s\ppq t$,
\item $N_{M}^{p+q-1}$ with $p+q > 2$ and $\min(s,t)=1$, $\max(s,t)=M$.
\end{enumerate}
\item An algebra of the form $\Gamma(p,q;r)$ (with $1 \ppq p \ppq q$) is stably equivalent of Morita
  type to an algebra of the
 form $\Lambda(p,q;s,t)$ if and only if they are isomorphic.
 This is only the case for $\Gamma(1,1;1)\cong \Lambda(1,1;2,2)$ and $\car K \neq 2$.
\item The algebras $\Gamma(p,q;r)$ and $\Gamma(p',q';r')$ (with $1 \ppq p\ppq q$ and $1 \ppq p'\ppq
  q'$) are stably  equivalent of Morita type if and only if $(p,q,r)=(p',q',r').$
\end{enumerate} \et

\begin{proof} It was proved by Rickard \cite{R} and Keller-Vossieck \cite{KV} that if two selfinjective
  $K$-algebras are derived equivalent, then they are stably equivalent of Morita type. Therefore to
  prove the result, we need only prove that the algebras listed in Theorem \ref{thm:derived equiv} are not
  stably equivalent of Morita type since they are all selfinjective.

We shall use the following invariants of stable equivalences of Morita type repeatedly:
\begin{itemize}
\item the dimension of $\HH^n(\Lambda)$ for $n\pgq 1$ for artin $K$-algebras (Xi, \cite[Theorem 4.2]{Xi});
\item the number of simple $\Lambda$-modules, where $\Lambda$ is an indecomposable selfinjective
  special biserial algebra (Pogorza\l y, \cite{Po});
\item the dimension of the centre $Z(\Lambda)\cong \HH^0(\Lambda)$ where $\Lambda$ is an indecomposable symmetric
  special biserial algebra (using a result of Liu, Zhou and Zimmermann, \cite[Corollary 1.2]{LZZ});
\item the absolute value of the Cartan determinant of $\Lambda$ (Xi, \cite[Proposition 5.1]{Xi}).
\end{itemize} All our algebras are indecomposable symmetric
  special biserial algebras. We now prove the theorem.
\begin{enumerate}[\bf (1)]
\item Assume first that $\Lambda(1,a;s,t)$ and $\Lambda(1,b;s',t')$, with $a\pgq1,$ $b\pgq1$, $2\ppq
  s\ppq t$ and $2\ppq s'\ppq t'$, are stably equivalent of Morita type. Then $a=b$ since the numbers
  of simples are equal and, using the dimensions of the centres, we have $s+t=s'+t'$. Moreover,
  their Cartan determinants are equal (they are positive) so that $st=s't'$. Finally, we have $(s,t)=(s',t')$ so that
  $\Lambda(1,a;s,t)=\Lambda(1,b;s',t')$.

Now assume that $N_M^a$ and $N_{M'}^b$,  with $a>1,$ $b>1$, $M\pgq 1$ and $M'\pgq 1$, are stably equivalent of Morita type. Then $a=b$ since the numbers
  of simples are equal and, using the Cartan determinant, we have $1+M+(a-1)M=1+M'+(a-1)M'$ so that
  $M=M'.$ Therefore $N_M^a=N_{M'}^b.$

Finally, if $\Lambda(1,a;s,t)$  and $N_{M}^b$,  with $a\pgq1,$ $b>1$, $2\ppq
  s\ppq t$ and $M\pgq 1$,  are stably equivalent of Morita type then, since $N_{M}^b$ is a symmetric
  Nakayama algebra, $\Lambda(1,a;s,t)$ must be a Brauer tree algebra by \cite{GR}, and hence derived
  equivalent to a symmetric Nakayama algebra by \cite{MH}, which we know is not the case. Therefore  $\Lambda(1,a;s,t)$  and $N_{M}^b$  are not stably equivalent of Morita type.

This concludes the proof of {\bf(1)}.
\item The proof is essentially the same as in Theorem \ref{thm:derived equiv} since almost all the
  derived invariants used there are invariants of stable equivalence of Morita type between
  indecomposable symmetric special biserial algebras, the exception being the K\"ulshammer ideal
  $T_1^\perp(A)$ of an algebra $A$. However, for symmetric algebras,
  the algebra $Z(A)/T_1^\perp(A)$ is a stable invariant of Morita type. Indeed, 
 let
  $Z^{st}(A)=\underline{\End}_{A^e}(A)$ be the
  stable centre of $A$ (the endomorphisms of $A$ in the stable $A^e$-module category)  and let  $Z^{pr}(A)=\ker({\End}_{A^e}(A)\rightarrow\underline{\End}_{A^e}(A))$ be
  the projective centre of $A.$ Then $Z^{st}(A)$ and $T_1^\perp(A)/Z^{pr}(A)$ are invariants of
  stable equivalences of Morita type for symmetric algebras (see \cite{LZZ,KLZ}), and  moreover
  $Z(A)/T_1^\perp(A)\cong Z^{st}(A)/(T_1^\perp(A)/Z^{pr}(A))$.

 Therefore the proof of Theorem
  \ref{thm:derived equiv} {\bf (2)} still holds for stable equivalences of Morita type.
\item The proof is the same as in Theorem \ref{thm:derived equiv} since the dimensions of the
  Hochschild cohomology groups in positive degrees are invariants of stable equivalences of Morita type.
\qedhere\end{enumerate}
\end{proof}

\br Recall from \cite{GR} that the Nakayama algebras are the distinct representatives of the
  stable equivalence classes of Brauer tree algebras. It follows from \cite{Po} and \cite{Liu} that any algebra
  which is stably equivalent of Morita type to one of the $\Gamma(p,q;r)$ or $\Lambda(p,q;s,t)$ is a 
  symmetric Brauer graph algebra. Moreover, the list of algebras given in Theorem \ref{thmSEMT} are
"normal forms" for derived equivalences of Brauer
graph algebras (a specific type of generalised star as in \cite{K} and  \cite[Theorem 5.7]{Ro}), and hence, since Brauer graph
algebras are selfinjective, for  stable equivalences of Morita type. Indeed, the algebra
$\Gamma(p,q;r)$ is the Brauer graph algebra associated to the graph  
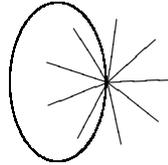
\begin{center}
$\xy /l2pc/:,
{\xypolygon9"P"{~={22}
~<{-}~>{}{}}};(0.89,0)*\ellipse(0.75,1.25){-};
 \endxy$\\\medskip 
\end{center}
with $p-1$ edges inside the loop and $q-1$ edges outside the loop and with multiplicity $r$ at the central vertex.
However, it is still an open question in general
whether two such normal forms are derived equivalent or not.
\er

\section*{Acknowledgements}

The first author would like to thank the  Universit\'e Jean Monnet (Saint-Etienne) and the Universit\'e Blaise Pascal (Clermont-Ferrand 2) for two months of invited professorships during which this collaboration took place.

\end{document}